\documentclass[11pt]{article}
\usepackage[colorlinks=true,
linkcolor=webgreen,
filecolor=webbrown,
citecolor=webgreen]{hyperref}

\definecolor{webgreen}{rgb}{0,.5,0}
\definecolor{webbrown}{rgb}{.6,0,0}

\usepackage{amssymb,psfig,graphics}

\newtheorem{theo}{Theorem}

\newtheorem{lem}{Lemma}

\newcommand{\Eqa}{\mbox{\footnotesize $t=1$}}
\newcommand{\Eqb}{\mbox{\footnotesize $t=-1$}}
\newcommand{\Eqc}{\mbox{\footnotesize $t=0$}}
\newcommand{\Eqd}{\mbox{\footnotesize $A=0$}}
\newcommand{\Eqe}{\mbox{\footnotesize $B=0$}}
\newcommand{\Eqf}{\mbox{\footnotesize $A=B$}}
\newcommand{\Eqg}{\mbox{\footnotesize $A=0$}}
\newcommand{\oaeq}{\mbox{\footnotesize ${\binom{1}{A}}$}}
\newcommand{\moaeq}{\mbox{\footnotesize $- {\binom{1}{A}}$}}
\newcommand{\ooxy}{\mbox{\footnotesize ${\binom{1}{X}} - {\binom{1}{Y}}$}}
\newcommand{\tob}{\mbox{\footnotesize $2 {\binom{1}{B}}$}}
\newcommand{\pooxy}{\mbox{\footnotesize ${\binom{1}{X}} + {\binom{1}{Y}}$}}
\newcommand{\mooxy}{\mbox{\footnotesize $- {\binom{1}{X}} - {\binom{1}{Y}}$}}
\newcommand{\treq}{\mbox{\footnotesize $tr \, \frac{b}{a} =0$}}
\newcommand{\vla}{v^{( \lambda )}}
\newcommand{\ta}{\tilde{a}}
\newcommand{\bone}{{\bf 1}}

\newcommand{\NR}{Nordstrom-Robinson~}
\newcommand{\Pre}{Preparata~}
\newcommand{\Ker}{Kerdock~}
\newcommand{\DG}{Delsarte-Goethals~}
\newcommand{\For}{Frobenius~}
\newcommand{\swe}{\mbox{swe}}
\newcommand{\cwe}{\mbox{cwe}}
\newcommand{\Lee}{\mbox{Lee}}
\newcommand{\Ham}{\mbox{Ham}}
\newcommand{\Cp}{{\sC^\perp}}
\newcommand{\PF}{\noindent{\em Proof}.~}
\newcommand{\ZZ}{{\mathbb Z}}
\newcommand{\sC}{{\mathcal C}}
\newcommand{\sG}{{\mathcal G}}
\newcommand{\sI}{{\mathcal I}}
\newcommand{\sP}{{\mathcal P}}
\newcommand{\sA}{{\mathcal A}}
\newcommand{\sK}{{\mathcal K}}
\newcommand{\sR}{{\mathcal R}}
\newcommand{\sE}{{\mathcal E}}
\newcommand{\sT}{{\mathcal T}}
\newcommand{\sH}{{\mathcal H}}
\newcommand{\la}{\lambda}
\newcommand{\af}{\alpha}
\def\binom#1#2{{#1}\choose{#2}}
\newcommand{\eqn}[1]{(\ref{#1})}
\newcommand{\hsp}{\hspace*{\parindent}}
\newcommand{\vsp}{\vspace{.1in}}
\newcommand{\bary}{\begin{eqnarray}}
\newcommand{\eary}{\end{eqnarray}}
\newcommand{\eeq}{\end{equation}}
\newcommand{\beql}[1]{\begin{equation}\label{#1}}
\newcommand{\bara}{\begin{eqnarray*}}
\newcommand{\eara}{\end{eqnarray*}}
\newcommand{\df}{\displaystyle\frac}

\newcommand{\dis}{\displaystyle}

\makeatletter
\def\section{\@startsection {section}{1}{\z@}{-3.5ex plus -1ex minus 
 -.2ex}{2.3ex plus .2ex}{\normalsize\bf}}
\def\subsection{\@startsection {subsection}{1}{\z@}{-3.5ex plus -1ex minus
 -.2ex}{2.3ex plus .2ex}{\normalsize\bf}}
\def\@sect#1#2#3#4#5#6[#7]#8{\ifnum #2>\c@secnumdepth
     \def\@svsec{}\else 
     \refstepcounter{#1}\edef\@svsec{\csname the#1\endcsname.\hskip .75em }\fi
     \@tempskipa #5\relax
      \ifdim \@tempskipa>\z@ 
        \begingroup #6\relax
          \@hangfrom{\hskip #3\relax\@svsec}{\interlinepenalty \@M #8\par}%
        \endgroup
       \csname #1mark\endcsname{#7}\addcontentsline
         {toc}{#1}{\ifnum #2>\c@secnumdepth \else
                      \protect\numberline{\csname the#1\endcsname}\fi
                    #7}\else
        \def\@svsechd{#6\hskip #3\@svsec #8\csname #1mark\endcsname
                      {#7}\addcontentsline
                           {toc}{#1}{\ifnum #2>\c@secnumdepth \else
                             \protect\numberline{\csname the#1\endcsname}\fi
                       #7}}\fi
     \@xsect{#5}}
\def\@begintheorem#1#2{\it \trivlist \item[\hskip \labelsep{\bf #1\ #2.}]}
\makeatother

\begin{document}
\begin{center}
{\large {\bf The $\ZZ_4$-Linearity of Kerdock, Preparata, Goethals and Related Codes\footnote{A different version of this paper appeared in: {\em IEEE Trans. Inform. Theory}, {\bf 40} (1994), 301--319.}}}\\
\vspace{\baselineskip}
{\em A. Roger Hammons, Jr.}** \\
\vspace{0.25\baselineskip}
Hughes Aircraft Company \\
Network Systems Division, Germantown, MD 20876 U.S.A. \\
\vspace{1\baselineskip}
{\em P. Vijay Kumar}** \\
\vspace{0.25\baselineskip}
Communication Science Institute, EE-Systems \\
University of Southern California, Los Angeles, CA 90089~U.S.A. \\
\footnotetext[7]{The work of A.~R. Hammons, Jr. and P.~V. Kumar was supported in part by the
National Science Foundation under Grant NCR-9016077 and by Hughes
Aircraft Company under its Ph.D. fellowship program.~~}
\vspace{\baselineskip}
{\em A.~R. Calderbank} and {\em N.~J.~A. Sloane} \\
\vspace{0.25\baselineskip}
Mathematical Sciences Research Center \\
AT\&T Bell Laboratories, Murray Hill, NJ 07974~U.S.A. \\
\vspace{\baselineskip}
{\em Patrick Sol\'{e}}$^{\S}$ \\
\vspace{0.25\baselineskip}
CNRS -- I3S, 250 rue A. Einstein, b\^{a}timent 4 \\
Sophia -- Antipolis, 06560 Valbonne, France \\
\footnotetext[4]{P. Sol\'{e} thanks the DIMACS Center and the IEEE for travel support.~~}
\vspace{1.5\baselineskip}
{\bf ABSTRACT}
\vspace{.5\baselineskip}
\end{center}
Certain notorious nonlinear binary codes contain more codewords than any
known linear code.
These include the codes constructed by \NR, Kerdock, Preparata,
Goethals, and \DG.
It is shown here that all these codes can be very simply constructed
as binary images under the Gray map of linear codes over $\ZZ_4$, the integers
$\bmod~4$ (although this requires a slight modification of the Preparata and Goethals codes).
The construction implies that all these binary codes are
distance invariant.
Duality in the $\ZZ_4$ domain implies that the binary images
have dual weight distributions.
The Kerdock and `Preparata' codes are duals over $\ZZ_4$ --- and the \NR
code is self-dual --- which
explains why their weight distributions are dual to each other.
The Kerdock and `Preparata' codes are $\ZZ_4$-analogues of first-order
Reed-Muller and extended Hamming codes, respectively.
All these codes are extended cyclic codes over $\ZZ_4$,
which greatly simplifies encoding and decoding.
An algebraic hard-decision decoding algorithm is given for the
`Preparata' code and a Hadamard-transform soft-decision
decoding algorithm for the Kerdock code.
Binary first- and second-order Reed-Muller codes are also linear
over $\ZZ_4$, but extended Hamming codes of length $n \ge 32$ and the Golay
code are not.
Using $\ZZ_4$-linearity, a new family of distance regular graphs are
constructed on the cosets of the `Preparata' code.
\clearpage
\thispagestyle{empty}
\setcounter{page}{1}
\begin{center}
{\large {\bf The $\ZZ_4$-Linearity of Kerdock, Preparata, Goethals and Related Codes}}\\
\vspace{\baselineskip}
{\em A. Roger Hammons, Jr.}** \\
\vspace{0.25\baselineskip}
Hughes Aircraft Company \\
Network Systems Division, Germantown, MD 20876 U.S.A. \\
\vspace{1\baselineskip}
{\em P. Vijay Kumar}** \\
\vspace{0.25\baselineskip}
Communication Science Institute, EE-Systems \\
University of Southern California, Los Angeles, CA 90089~U.S.A. \\
\footnotetext[7]{The work of A.~R. Hammons, Jr. and P.~V. Kumar was supported in part by the
National Science Foundation under Grant NCR-9016077 and by Hughes
Aircraft Company under its Ph.D. fellowship program.~~}
\vspace{\baselineskip}
{\em A.~R. Calderbank} and {\em N.~J.~A. Sloane} \\
\vspace{0.25\baselineskip}
Mathematical Sciences Research Center \\
AT\&T Bell Laboratories, Murray Hill, NJ 07974~U.S.A. \\
\vspace{\baselineskip}
{\em Patrick Sol\'{e}}$^{\S}$ \\
\vspace{0.25\baselineskip}
CNRS -- I3S, 250 rue A. Einstein, b\^{a}timent 4 \\
Sophia -- Antipolis, 06560 Valbonne, France \\
\footnotetext[4]{P. Sol\'{e} thanks the DIMACS Center and the IEEE for travel support.~~}
\vspace{1.5\baselineskip}
\end{center}
\setlength{\baselineskip}{1.5\baselineskip}
\section{Introduction}
\hspace*{\parindent}
Several notorious families of nonlinear codes have more codewords
than any comparable linear code presently known.
These are the \NR, Kerdock, \Pre, Goethals and
\DG codes
\cite{BT92},
\cite{DG75},
\cite{Go74},
\cite{Go76},
\cite{Ke72},
\cite{MS77},
\cite{NR67},
\cite{Pr68}.
Besides their excellent error-correcting capabilities these codes
are remarkable because the Kerdock and Preparata codes are
`formal duals', in the sense that although these codes are nonlinear,
the weight distribution of one is the MacWilliams
transform of the weight distribution of the other \cite[Chap.~15]{MS77}.
The main unsolved question concerning these codes has always been whether
they are duals in some more algebraic sense.
Many authors have investigated these codes, and have found that
(except for the Nordstrom-Robinson code)
they are not unique, and indeed that large numbers of codes exist
with the same weight distributions
\cite{BLW83}, \cite{Ca89},
\cite{Ka82},
\cite{Ka82a},
\cite{Ka83},
\cite{Li83}.
Kantor \cite{Ka83} declares that the ``apparent relationship between these
[families of codes] is merely a coincidence.''

Although this may be true for many versions of these codes,
we will show that, when properly defined,
Kerdock and Preparata codes are {\em linear} over
$\ZZ_4$ (the integers $\bmod~4$), and that as $\ZZ_4$-codes
they {\em are} duals.
They are in fact just extended cyclic codes over $\ZZ_4$.

The version of the Kerdock code that we use is the standard one,
while our version of the Preparata code differs from the standard
one in that it is not a subcode of the extended Hamming code but of a
nonlinear code with the same weight distribution as the extended
Hamming code.
Our `Preparata' code has the same weight distribution as Preparata's version,
and has a similar construction in terms of finite field
transforms.
In our version, the Kerdock and `Preparata' codes are
$\ZZ_4$-analogues of first-order Reed-Muller and extended Hamming codes,
respectively.
Since the new construction is so simple, we propose that this is the
`correct' way to define these codes.

The situation may be compared with that for Hamming codes.
It is known that there are many binary codes with the same weight distribution
as the Hamming code --- all are perfect single-error correcting codes,
but one is distinguished by being linear
(see \cite{Va62}, \cite{Ph84}, \cite{Phxx} and also \S5.4).
Similarly, there are many binary codes with the same weight distributions as the
Kerdock and \Pre codes;
one pair is distinguished by being the images of a dual pair of 
linear extended-cyclic codes over $\ZZ_4$.
It happens that \Ker picked out the distinguished code,
although \Pre did not.

Kerdock and Preparata codes exist for all lengths $n=4^m \ge 16$.
At length 16 they coincide,
giving the Nordstrom-Robinson code \cite{NR67},
\cite{Sn73}, \cite{Go77}.
The $\ZZ_4$ version of the Nordstrom-Robinson code turns out to be the
`octacode' \cite{CS92}, \cite{CS93}, a self-dual code of length 8 over $\ZZ_4$ that is used when the Leech lattice is constructed from eight copies of
the face-centered cubic lattice.

The very good nonlinear binary codes of minimal distance 8 discovered by
Goethals \cite{Go74}, \cite{Go76}, and the high minimal distance codes of
Delsarte and Goethals \cite{DG75},
also have a simple description as extended cyclic codes over $\ZZ_4$,
although our `Goethals' code differs slightly from Goethals' original
construction.

The decoding of all these codes is greatly simplified by working in the $\ZZ_4$-domain,
where they are linear and it is meaningful to speak of syndromes.
Decoding the \NR and `Preparata' codes is especially simple.

These discoveries came about in the following way.
Recently, a family of nearly optimal four-phase sequences of period
$2^{2r+1} -1$, with alphabet $\{1, i, -1, -i \}$, $i = \sqrt{-1}$,
was discovered by Sol\'{e} \cite{So89} and later independently by Bozta\c{s},
Hammons and Kumar \cite{Bo90}, \cite{BHK92}.
By replacing each element $i^a$ by its exponent $a \in \{0,1,2,3\}$,
this family may be viewed as a linear code over
$\ZZ_4$.
Since the family has low correlation values, it also possesses a large
minimal Euclidean distance and thus has the potential for excellent
error-correcting capability.

When studying these four-phase sequences, Hammons and Kumar
and later independently Calderbank, Sloane and Sol\'{e} noticed the
striking resemblance between the 2-adic (i.e. base~2)
expansions of the quaternary codewords and the standard construction
of the \Ker codes.
The reader can see this for himself by comparing the formulae on
page 1107 of \cite{BHK92} (the common starting point for the two
independent discoveries) and page 458 of \cite{MS77}.

Both teams then realized that the \Ker code is simply the image of the quaternary code (when extended by an zero-sum check symbol) under the Gray map
defined below (see \eqn{E1}).
This was a logical step to pursue since the Gray map translates a
quaternary code with high minimal Lee or Euclidean distance into a binary code of
twice the length with high minimal Hamming distance.

The discovery that the quaternary dual gives a code which is the `correct' definition of the `Preparata' code
followed almost immediately.

The two teams worked independently until the middle of November 1992,
when, discovering the considerable overlap between their work,
they decided to join forces.
The discoveries about the \Ker and \Pre codes are
in a paper \cite{xx} presented by Hammons and Kumar at the International Symposium on
Information Theory (San Antonio, January 1993, but submitted in June
1992),
in Hammons' dissertation \cite{Ha92},
and in a manuscript \cite{xy}  (now replaced by the present paper) submitted in
early November 1992 to these {\em Transactions}.
Hammons and Kumar realized in June 1992 that the $\ZZ_4$ \Ker and
`Preparata' codes could be generalized to give the quaternary Reed-Muller codes $QRM (r,m)$ of Section~5.4.

In late October 1992,
Calderbank, Sloane and Sol\'{e} submitted a research
announcement (now replaced by \cite{CHKSS}) to the
{\em Bulletin of the American Mathematical Society}, also
containing the discoveries about the \Ker and \Pre codes,
as well as results (Sections~2.6 to 2.8) about the existence of quaternary versions of Reed-Muller,
Golay and Hamming codes.
They discovered the quaternary versions of the Goethals and \DG codes
in early November.

The present paper is a compositum of all our results.

The discovery that the \NR code is a quaternary version of the octacode was
made by Forney, Sloane and Trott in early October 1992,
and is described in \cite{FST93}.
(It was already known to Hammons and Kumar in June 1992 that the Nordstrom-Robinson code was linear over $\ZZ_4$, but they had not made the identification
with the octacode.)

It can be shown that the binary nonlinear single-error-correcting codes
found by Best \cite{Be80}, Julin \cite{Ju65}, Sloane and Whitehead \cite{SW70}
and others can also be more simply described as codes over
$\ZZ_4$ (although here the corresponding $\ZZ_4$-codes are
nonlinear).
This will be described elsewhere \cite{CS93b}.
Large sequence families for code-division multiple-access (CDMA)
that are supersets of the near optimum four-phase sequence
families described above and which are related to the
Delsarte-Goethals codes are investigated in \cite{xz}.

The paper is arranged as follows.
Section~II discusses linear codes over $\ZZ_4$, their duals,
and their images as binary codes under the Gray map.
Necessary and sufficient conditions are given for a binary code to be
the image of a linear code over $\ZZ_4$.
Reed-Muller codes of length $2^m$ and orders $0,1,2, m-1 ,m$ satisfy
these conditions, but extended Hamming codes and the Golay code do not.
Cyclic codes over $\ZZ_4$ are studied by means of Galois rings $GR (4^m)$
rather than the Galois fields $GF(2^m)$ used to analyze binary cyclic
codes,
and
Section~III is devoted to these rings.

In Section~IV we show that \Ker codes are extended cyclic codes
over $\ZZ_4$, and in fact are simply $\ZZ_4$-analogues of first-order
Reed-Muller
codes (see the generator matrix \eqn{E3.12} and also \S5.4).
The \NR code is discussed in \S4.5.
Subsequent
subsections give the weight distribution of the \Ker
codes and a soft-decision decoding algorithm for them.

In Section~V
we show that the binary images of the quaternary duals of the \Ker codes are
Preparata-like codes, having essentially the same properties as Preparata's
original codes.
Theorem~\ref{PD}, however, shows that the `Preparata' codes are
strictly different from the original construction.
\S5.2 provides a finite field transform characterization of the `Preparata' codes and compares them with the original codes.
The `Preparata' codes have a very simple decoding algorithm (\S5.3).
(This is considerably simpler than any previous decoding
algorithm --- compare \cite{Bo88}.)
Section~5.4 defines a family of quaternary Reed-Muller codes $QRM (m,r)$ which generalizes the quaternary Kerdock and `Preparata' codes.
The final subsections are concerned with the automorphism groups of
these codes (\S5.5),
and a new family of distance regular graphs defined on the cosets of
the `Preparata' code (\S5.6).

In Section~VI we show that the binary nonlinear
\DG codes \cite{DG75}
are also
extended cyclic codes over $\ZZ_4$,
and that their $\ZZ_4$-duals have essentially the same properties
as the Goethals codes \cite{Go74}, \cite{Go76} and the
`Goethals-Delsarte' codes of Hergert \cite{He90}.

\vspace{.2in}
\noindent{\bf Postscript.}
After this paper was completed, V.~I. Levenshtein
drew our attention to an article by
Nechaev \cite{Nec}.
In this article Nechaev considers the quaternary sequences
$\{c_t \}$ given (in the notation of the present paper) by
$$c_t = (-1)^t \{ T( \la \xi^t ) + \delta \} ~,$$
$0 \le t \le 2^{m+1} -3$, $\la \in R$, $\delta \in \ZZ_4$,
and their 2-adic expansions $c_t= a_t + 2b_t$, where
$a_t$, $b_t \in \{0,1\}$.
The principal result of \cite{Nec} shows
that the set of $\{b_t\}$ is a nonlinear
binary cyclic code which is equivalent to the binary
Kerdock code punctured in two coordinates.
However, \cite{Nec} makes no mention of the fundamental isometry of Eq.~\eqn{E1}, nor of Preparata codes and the sense in which they are duals
of Kerdock codes.
\section{Quaternary and related binary codes}
\subsection{Quaternary codes.~~}
By a {\em quaternary code} $\sC$ of length $n$
we shall mean a linear block code over $\ZZ_4$,
i.e. an additive subgroup of $\ZZ_4^n$.
Such codes have been studied recently both in
connection with the construction of sequences with low correlation
(\cite{Bo90}, \cite{BHK92}, \cite{So89}, \cite{US91})
and in a variety of other contexts (see \cite{CS93} and the
references contained therein).

We define an inner product on $\ZZ_4^n$ by $a \cdot b = a_1 b_1 + \cdots + a_n b_n$ $( \bmod~4)$, and then the notions of {\em dual} code $(\sC^\perp )$,
{\em self-orthogonal} code $( \sC \subseteq \sC^\perp )$ and
{\em self-dual} code
$( \sC = \sC^\perp )$ are defined in the standard way
(cf. \cite{Kl87}, \cite{MS77}).
For many applications there is no need to distinguish between $+1$
components of codewords and $-1$ components, and so we say that two codes are
{\em equivalent} if one can be obtained from the other by
permuting the coordinates and (if necessary)
changing
the signs of certain coordinates.
Codes differing by only a permutation of coordinates are called
{\em permutation-equivalent}.
The automorphism group $\mbox{Aut} ( \sC )$ of $\sC$ consists of all
permutations and sign-changes of the coordinates that preserve the set of
codewords.

Any code is permutation-equivalent to a code $\sC$ with generator
matrix of the form
\beql{E2.1}
G = \left[ \begin{array}{ccc}
I_{k_1} & A & B \\
0 & 2I_{k_2} & 2C
\end{array}
\right] ~,
\eeq
where $A$ and $C$ are $\ZZ_2$-matrices and $B$ is a $\ZZ_4$-matrix.
The code is then an elementary abelian group of type $4^{k_1} 2^{k_2}$,
containing $2^{2k_1 + k_2}$ codewords.
We shall indicate this by
saying that $\sC$ has type $4^{k_1} 2^{k_2}$, or simply that
$| \sC | = 4^{k_1} 2^{k_2}$.

Eq.~\eqn{E2.1} illustrates a difference in point of view between ring
theory and coding theory.
Quaternary codes are $\ZZ_4$-modules.
A ring theorist would point out, correctly, that a quaternary code is not
in general a free module \cite{Hu74},
and so need not have a basis.
Although this is true, \eqn{E2.1}
is a perfectly good generator matrix.
Encoding is carried out by writing the information symbols
in the form $u= u_1 \cdots u_{k_1} u_{k_1 +1} \cdots u_{k_1 + k_2}$, where
$u_i \in \ZZ_4$
if $1 \le i \le k_1$, $u_i \in \ZZ_2$ if
$k_1 +1 \le i \le k_1 + k_2$, and mapping $u$ to the codeword
$uG$.
The code
$\sC$ is a free $\ZZ_4$-module if and only if $k_2 =0$.

If $\sC$ has generator matrix \eqn{E2.1}, the
dual code $\sC^\perp$ has generator matrix
\beql{E2.2}
\left[
\begin{array}{ccc}
-B^{tr} -C^{tr}A^{tr} & C^{tr} & I_{n-k_1 -k_2} \\
~~ \\ [-.1in]
2A^{tr} & 2I_{k_2} & 0
\end{array}
\right]
\eeq
and type $4^{n-k_1 - k_2} 2^{k_2}$.
\subsection{Weight enumerators.~~}
Several weight enumerators are associated with a quaternary code $\sC$.
The {\em complete weight enumerator} (or c.w.e.) of $\sC$ is
\beql{E2.3}
\cwe_\sC (W,X,Y,Z) =
\sum_{a \in \sC} W^{n_0 (a)} X^{n_1 (a)} Y^{n_2 (a)}
Z^{n_3 (a)} ~,
\eeq
where $n_j (a)$ is the number of components of $a$ that are congruent to $j$
($\bmod~4$) (cf.~\cite{Kl87}, \cite[p.~141]{MS77}).
Permutation-equivalent codes have the same c.w.e., but equivalent codes may have distinct c.w.e.'s.
The appropriate weight enumerator for an equivalence class of codes is
the {\em symmetrized weight enumerator} (or s.w.e.), obtained by
identifying $X$ and $Z$ in \eqn{E2.3}:
\beql{E2.4}
\swe_\sC (W,X,Y) = \cwe_\sC (W,X,Y,X) ~.
\eeq

The Lee weights of $0,1,2,3 \in \ZZ_4$ are $0,1,2,1$ respectively,
and the Lee weight $wt_L (a)$ of $a \in \ZZ_4^n$ is
the rational sum of the Lee weights of its components.
This weight function defines a distance $d_L (~,~)$ on $\ZZ_4^N$ called
the {\em Lee metric}.
The {\em Lee weight enumerator} of $\sC$ is
\begin{eqnarray}
\label{E2.12}
\Lee_\sC (W,X) & = &
\sum_{a \in \sC} W^{2n- wt_L (a)}
X^{wt_L (a)} \nonumber \\
& = & \swe_\sC ( W^2 , WX , X^2 ) ~,
\end{eqnarray}
a homogeneous polynomial of degree $2n$.
Finally, the {\em Hamming weight enumerator} of $\sC$,
less useful than the others, is
\beql{E2.13}
\mbox{Ham}_\sC (W,X) =
\swe_\sC (W,X,X) ~.
\eeq

We then have the following analogues of the MacWilliams identity,
giving the weight enumerators for the dual code $\sC^\perp$
(\cite{Kl87}, \cite{Kl89}, \cite{CS93}):
\beql{E2.5}
\cwe_{\sC^\perp} (W,X,Y,Z) = \frac{1}{| \sC |}
\cwe_{\sC} (W+X+Y+Z, W+ iX - Y - iZ ,
W-X+Y-Z, W- iX -Y+iZ ) ~,
\eeq
\beql{E2.6}
\swe_{\sC^\perp} (W,X,Y) = \frac{1}{| \sC|} \swe_{\sC}
(W + 2X + Y , W-Y , W-2X+Y) ~,
\eeq
\beql{E2.7}
\Lee_{\Cp} (W,X) = \frac{1}{| \sC|} \Lee_{\sC}
(W+X, W-X) ~,
\eeq
\beql{E2.8}
\mbox{Ham}_{\Cp} (W,X) =
\frac{1}{| \sC|}
\mbox{Ham}_{\sC} (W+3X , W-X) ~.
\eeq
There are also several analogues of Gleason's theorem, giving bases for the
weight enumerators of self-dual codes --- see \cite{Kl87}, \cite{CS93}.

\subsection{Associated complex-valued sequences.~~}
We may associate to every $\ZZ_4$-valued vector
$a = (a_1 , \ldots , a_n )$ an equivalent complex roots-of-unity sequence
$s=i^a = (i^{a_1} , \ldots , i^{a_n} )$, where $i= \sqrt{-1}$.
Then, given a set $\sC$ of quaternary vectors,
we let
$$\Omega ( \sC ) = \{ i^a :~ a \in \sC \}$$
denote the corresponding set of complex sequences.
When $\sC$ is regarded as a set of CDMA
signature sequences, its effectiveness depends on the complex
correlations (or Hermitian inner products) of the sequences
in $\Omega ( \sC )$.
When $\sC$ is regarded as a code, its error-correcting
capability depends on the Euclidean distance
properties of $\Omega ( \sC )$.
If $a,b$ are quaternary vectors with associated vectors
$s=i^a$, $t=i^b$, then
\begin{eqnarray}
\label{H2.1}
\|s-t \|^2 & = & \|s \|^2 + \|t \|^2 - 2 ~ \mbox{Re}~
\{ s^H t \}  \nonumber \\
& = & 2n - 2 ~ \mbox{Re} ~ \{ \zeta (a-b) \} ~,
\end{eqnarray}
where $^H$ denotes the Hermitian inner product, and
\beql{H2.2}
\zeta (a-b) = \sum_{r=1}^n i^{a_r - b_r}
\eeq
is the {\em complex correlation} of $a$ and $b$.
Note that $\zeta$ depends only on the difference $a-b$.
By \eqn{H2.1}, if the nontrivial correlations of $\Omega ( \sC )$ are low in magnitude, then the set
possesses large minimal Euclidean distance.
We also see that
\beql{E2.9}
\| s-t \|^2 = 2d_L (a,b) ~.
\eeq
\subsection{Binary codes associated with quaternary codes; the Gray map.~~}
In communication systems employing quadrature phase-shift keying (QPSK),
the preferred assignment of two information bits to the four possible
phases is the one shown in Fig.~1,
in which adjacent phases differ by only one binary digit.
This mapping is called {\em Gray encoding} and has the advantage
that, when a quaternary codeword is transmitted across an additive
white Gaussian noise channel, the errors most likely
to occur are those causing a single erroneously decoded information bit.
\begin{figure}[htb]
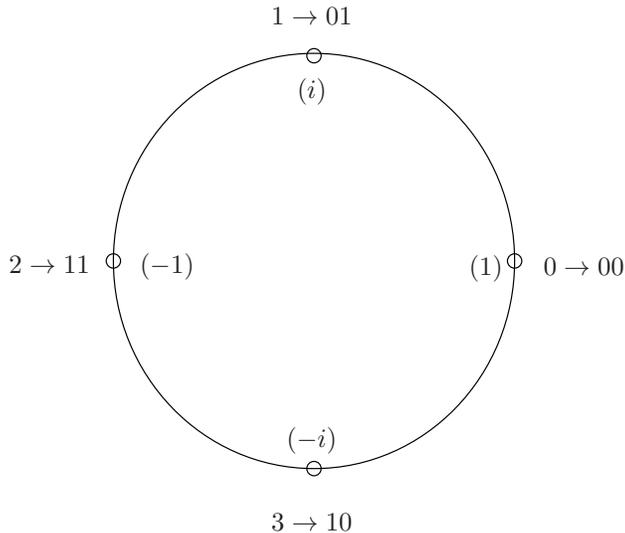

\begin{center}
\input SD1.pstex_t
\end{center}
\caption{Gray encoding of quaternary symbols and QPSK phases.~~}
\end{figure}

Formally, we define three maps from $\ZZ_4$ to $\ZZ_2$ by
$$
\begin{array}{cccc}
c & \af (c) & \beta (c) & \gamma (c) \\
~~  \\ [-.09in]
0 & 0 & 0 & 0 \\
1 & 1 & 0 & 1 \\
2 & 0 & 1 & 1 \\
3 & 1 & 1 & 0
\end{array}
$$
and extend them in the obvious way to maps from $\ZZ_4^n$
to $\ZZ_2^n$.
The 2-adic expansion of $c \in \ZZ_4$ is
\beql{E0}
c = \af (c) + 2 \beta (c) ~.
\eeq
Note that
$\af (c) + \beta (c) + \gamma (c) =0$ for all $c \in \ZZ_4$.
We construct binary codes from quaternary codes using the
{\em Gray map} $\phi: \ZZ_4^n \to \ZZ_2^{2n}$ given by
\beql{E1}
\phi (c) = ( \beta (c) , \gamma (c)), ~~~
c \in \ZZ_4^n ~.
\eeq
When we speak of the binary image of a quaternary code $\sC$,
we will always mean its image $C= \phi ( \sC )$ under the Gray map.
We use script letters for quaternary codes, with the corresponding
Latin letters for their binary images.

$C$ is in general a nonlinear binary code of length $2n$.
If $C$ {\em is} linear, and $\sC$ is defined by
\eqn{E2.1},
then $C$ has generator matrix
\beql{E32}
\left[
\begin{array}{cccccc}
I_{k_1} & A & \af (B) & I_{k_1} & A & \af (B) \\
~~ \\ [-.09in]
0 & I_{k_2} & C & 0 & I_{k_2} & C \\
~~ \\ [-.09in]
0 & 0 & \beta (B) & I_{k_1} & A & \gamma (B)
\end{array}
\right] ~.
\eeq

We say that a binary code $C$ is {\em $\ZZ_4$-linear} if its
coordinates can be arranged so that it is the image under the Gray map
$\phi$ of a quaternary code $\sC$.

The crucial property of the Gray map is that it preserves distances.

\begin{theo}
\label{Iso}
$\phi$ is a distance-preserving map or {\em isometry} from
$$
(\ZZ_4^n , ~ \mbox{Lee distance}) ~~\mbox{to}~~
( \ZZ_2^{2n} , ~ \mbox{Hamming distance}) ~.
$$
\end{theo}

\noindent{\em Proof}.~
It is easy to see from the definitions (and Fig.~1) that
\begin{eqnarray}
\label{E2.10}
wt ( \phi (a)) & = & wt_L (a) , ~~~~ a \in \ZZ_4^n ~, \\
\label{E2.11}
d( \phi (a) , \phi (b)) & = &
d_L (a,b) , ~~~~a,b \in \ZZ_4^n ~,
\end{eqnarray}
where $wt(~)$ and $d(~,~)$ are the usual Hamming weight and
distance functions for binary vectors.~~~ \hfill $\blacksquare$

From \eqn{E2.9}, \eqn{E2.11},
the Hamming distance between the binary images
$\phi (a)$ and $\phi (b)$ is proportional to the squared Euclidean distance between
the complex sequences $i^a$ and $i^b$.

Two other binary codes $C^{(1)}$, $C^{(2)}$ are
canonically associated with a quaternary code $\sC$.
These are the linear codes defined by
\begin{eqnarray}
\label{E33}
C^{(1)} & = & \{ \af(c) : ~ c \in \sC \} ~, \\
\label{E34}
C^{(2)} & = & \{ \beta (c) : ~ c \in \sC , ~ \af (c) =0 \} ~.
\end{eqnarray}
If $\sC$ has generator matrix \eqn{E2.1}, then $C^{(1)}$ is an
$[n,k_1]$ code with generator matrix
\beql{E35}
[I_{k_1} ~~A~~ \af (B) ] ~,
\eeq
while $C^{(2)} \supseteq C^{(1)}$ is an $[n,k_1 + k_2 ]$ code with generator matrix
\beql{E36}
\left[
\begin{array}{ccc}
I_{k_1} & A & \af (B) \\
~~ \\ [-.09in]
0 & I_{k_2} & C
\end{array}
\right]
\eeq
--- compare \eqn{E32}.
It is shown in \cite{CS93} that given any binary codes
$C'$, $C''$ of length $n$ with $C'' \supseteq C'$,
there is a quaternary code $\sC$ with $C^{(1)} =C'$,
$C^{(2)} = C''$.
\subsection{Weight and distance properties.~~}
Since in general $C = \phi ( \sC )$ is not linear,
it need not have a dual.
We define its $\ZZ_4$-{\em dual} to be
$C_\perp = \phi ( \sC^\perp )$, as in the diagram
\beql{E2.17}
\begin{array}{r@{~}ccc@{~}c@{~}c}
~ & \sC & \stackrel{\phi}{\longrightarrow} & C & = & \phi ( \sC ) \\
~~ \\ [-.1in]
\mbox{dual} & \Biggl \downarrow \\
~~ \\ [-.1in]
~ &\sC^\perp & \stackrel{\phi}{\longrightarrow} & C_\perp & = & \phi ( \Cp ) ~.
\end{array}
\eeq
Note that one cannot add an arrow marked `dual' on the right side to produce
a commuting diagram.

In this section we discuss the weight and distance properties
of $C$ and $C_\perp$.
The principal results to be derived here are the following:
\begin{itemize}
\item[(1)]
$C$ and $C_\perp$ are distance invariant.
\item[(2)]
The weight distributions of $C$ and $C_\perp$ are MacWilliams
transforms of one another.
\end{itemize}

A binary code $C$ is said to be {\em distance invariant}
\cite[p.~40]{MS77} if the Hamming weight distributions
of its translates $u+C$ are the same for all $u \in C$.

\begin{theo}
\label{tH2.1}
If $\sC$ is a $($linear$)$ quaternary code, then its binary Gray representation $C= \phi ( \sC )$ is distance invariant.
\end{theo}

\noindent{\em Proof}.~
$\sC$ is distance invariant (with respect to Lee distance) because it
is linear, and the result then follows from Theorem~\ref{Iso}. \hfill $\blacksquare$

For a distance invariant code $C$ of length $n$, the
(Hamming) weight enumerator
$$\mbox{Ham}_C (W,Z) = \sum_{c' \in C} W^{n-d(c',c)} X^{d(c',c)}$$
is independent of $c \in C$.
If $C= \phi ( \sC )$, it follows from Theorem~\ref{Iso} and \eqn{E2.12} that
\beql{E2.14}
\mbox{Ham}_C (W,X) = \Lee_{\sC} (W,X) =
\swe_{\sC} (W^2 , WX ,X^2 ) ~.
\eeq

\begin{theo}
\label{H2.5}
If $\sC$ and $\Cp$ are dual quaternary codes, then the weight distributions
of the binary codes $C = \phi ( \sC )$ and $C_\perp = \phi ( \sC^\perp )$
are related by the binary MacWilliams transform.
\end{theo}

\noindent{\em Proof}.~
From \eqn{E2.14}, \eqn{E2.7} we have
\begin{eqnarray*}
\mbox{Ham}_{C_\perp} (W,X) & = &
\Lee_{\Cp} (W,X) \\
& = & \frac{1}{| \sC |} \Lee_\sC (W+X,~W-X) \\
& = & \frac{1}{| \sC |} \Ham_C (W+X ,~W-X) ~.
\end{eqnarray*}
as required. \hfill $\blacksquare$
\subsection{Existence and linearity conditions.~~}
We now give necessary and sufficient conditions for a binary code to be
$\ZZ_4$-linear, and for the binary image of a quaternary code to be
a linear code.
The reader who is primarily interested in \Ker and \Pre codes
should skip to Section~III.

Since $\phi (-c) = ( \gamma(c) , \beta (c))$,
it follows that if $C$ is $\ZZ_4$-linear then $C$ is fixed under
the `swap' map $\sigma$ that interchanges the left and right
halves of each codeword:
\beql{E1b}
\sigma : (u_1 ~ u_2 \cdots u_n ~ u_{n+1} \cdots u_{2n} ) \mapsto
(u_{n+1} \cdots u_{2n} ~ u_1 ~ u_2 \cdots u_n ) ~.
\eeq
In other words $\sigma$ applies the permutation
\beql{E1c}
(1,n+1) (2, n+2) \cdots (n,2n)
\eeq
to the coordinates.
This is a fixed-point-free involution in the automorphism group of $C$.

\begin{theo}
\label{L2}
A binary, not necessarily linear, code $C$ of even length is $\ZZ_4$-linear
if and only if its coordinates can be arranged so that
\beql{E4}
u,v \in C \Rightarrow u+v+ (u+ \sigma(u)) \ast (v+ \sigma (v)) \in C ~,
\eeq
where $\sigma$ is the swap map that interchanges the left and right
halves of a vector, and $\ast$ denotes the componentwise
product of two vectors.
\end{theo}

\noindent{\em Proof}.~
This is an immediate consequence of the easily-verified identity
\beql{E5}
\phi (a+b) = \phi (a) + \phi (b) + ( \phi (a) + \sigma ( \phi (a))) \ast
( \phi (b) + \sigma ( \phi (b))) ~,
\eeq
for all $a,b \in \ZZ_4^n$. \hfill $\blacksquare$
\begin{theo}
\label{L1}
The binary image $\phi ( \sC )$ of a quaternary linear code
$\sC$ is linear if and only if
\beql{E2}
a,b \in \sC \Rightarrow 2 \af (a) \ast \af (b) \in \sC ~
\eeq
\end{theo}

\noindent{\em Proof}.~
This is an immediate consequence of the identity
\beql{E3}
\phi (a) + \phi (b) + \phi (a+b) = \phi (2 \af (a) \ast \af (b))
\eeq
for all $a,b \in \ZZ_4^n$. \hfill $\blacksquare$
\begin{theo}
\label{L3}
A binary linear code $C$ of even length is $\ZZ_4$-linear if and only if
its coordinates can be permuted so that
\beql{E6}
u,v \in C \Rightarrow (u+ \sigma (u)) \ast (v+ \sigma (v)) \in C ~,
\eeq
where $\sigma$ is as in Theorem~\ref{L2}.
\end{theo}

\noindent{\em Proof}.~
This is also a consequence of \eqn{E5}. \hfill $\blacksquare$

Conditions \eqn{E2}, \eqn{E4} and \eqn{E6} are very restrictive,
and (we are now speaking informally)
imply that most binary codes are not $\ZZ_4$-linear.
\subsection{Reed-Muller and Hamming codes}

\begin{theo}
\label{RM}
The $r$th order binary Reed-Muller code $RM (r,m)$ of length $n=2^m$,
$m \ge 1$, is $\ZZ_4$-linear for $r=0,1,2, m-1$ and $m$.
\end{theo}

\noindent{\em Proof}.~
We leave to the reader the straightforward verification that
$RM(r,m)$ is the image under $\phi$ of the quaternary code
$ZRM (r,m-1)$ (say) of length $2^{m-1}$ generated by $RM(r-1, m-1)$ and
$2~RM (r,m-1)$, for
$r=0,1,2,m-1,m$
(with the convention that $RM (-1, m-1) = RM (m,m-1) =0$). \hfill $\blacksquare$

Let $(v_1 , \ldots , v_{m-1} )$ range over $\ZZ_2^{m-1}$,
so that $RM(r,m-1)$ is generated (in the usual way, as a binary code)
by the vectors corresponding to monomials in the Boolean functions
$v_i$ of degree $\le r$ \cite[Chap.~13]{MS77}.
Then $RM(1,m)$ is the binary image of the quaternary code
$ZRM (1,m-1)$ generated by the vectors
corresponding to $1, 2v_1 , \ldots , 2v_{m-1}$, and
$RM(2,m)$ is the image of the quaternary code $ZRM (2,m-1)$ generated by
$1, v_1 , \ldots , v_{m-1}$, $2v_1 v_2$, $2v_1 v_3 , \ldots , 2v_{m-2} v_{m-1}$.

For example, the $[16,5,8]$ code $RM (1,4)$ and the $[16,11,4]$ code
$RM(2,4)$ are the binary images of the quaternary codes
with generator matrices
\beql{E21}
\begin{array}{c}
ZRM (1,3) \\
\left[ \begin{array}{c}
11111111 \\
00002222 \\
00220022 \\
02020202
\end{array}
\right]
\end{array}
\begin{array}{l}
~~~ \\
1 \\
2v_1 \\
2v_2 \\
2v_3
\end{array} ~,~~~
\begin{array}{c}
ZRM (2,3) \\
\left[
\begin{array}{c}
11111111 \\
00001111 \\
00110011 \\
01010101 \\
00000022 \\
00000202 \\
00020002
\end{array}
\right]
\end{array}
\begin{array}{l}
~ \\
1 \\
v_1 \\
v_2 \\
v_3 \\
2v_1 v_2 \\
2v_1 v_3 \\
2v_2 v_3
\end{array} ~.
\eeq

In Eq.~\eqn{E6}, if $u,v$ are represented by Boolean functions of degree $r$,
and $(u+ \sigma (u)) \ast (v+ \sigma (v)) \neq 0$, then
$(u+ \sigma (u)) \ast (v+ \sigma (v))$ is a Boolean function of degree $2r-2$.
So an $r$th order Reed-Muller code with $r \le m/2$ satisfies
\eqn{E6} provided $r \le 2$ (which gives an alternative proof of part of
Theorem~\ref{RM}), but we conjecture that it does not
satisfy \eqn{E6} if $3 \le r \le m-2$.
In other words we conjecture that if $C$ is a binary Reed-Muller code
$RM (r,m)$ with $3 \le r \le m-2$, then there is no permutation of the
coordinates of $C$ such that the permuted code is equal to $\phi ( \sC )$ for
some quaternary code $\sC$.
However, we have found a proof of this only for $(m-2)$nd order RM codes.

\begin{theo}
\label{XH}
The binary code $RM(m-2,m)$, i.e. the extended Hamming code of length
$n=2^m$, is not $\ZZ_4$-linear for $m \ge 5$.
\end{theo}

\noindent{\em Proof}.~
Suppose $H$ is a $[2^m , 2^m - m-1 , 4]$ extended Hamming code with its coordinates arranged so that
$H= \phi ( \sH )$ for some quaternary code $\sH$.
We will obtain a contradiction for $m \ge 5$.
The codewords of weight 4 in $H$ form a Steiner system
$S(3,4,2^m)$ \cite[p.~63]{MS77}.
From this it follows without difficulty that
\beql{E30}
\begin{array}{l}
\mbox{for $m \ge 4$, $H$ contains codewords of} \\
\mbox{weight 4 that meet in just one coordinate.}
\end{array}
\eeq
Let $F$ be the subcode of $H$ fixed under the swap map $\sigma$ of
\eqn{E1b}, and let $\psi$ be the homomorphism $H \to F$
given by $\psi (x) = x+ \sigma (x)$.
Then $\mbox{im} \, \psi \subseteq \ker \psi =F$.
Since $\dim \ker \psi \le 2^{m-1} -1$, $\dim ~ \mbox{im}~ \psi \ge 2^{m-1} -m$.
Let $E$ consist of the right-hand halves of the codewords in $\mbox{im}~\psi$.
Then $E$ is a $[2^{m-1} , \ge 2^{m-1} -m , 2]$ code,
containing say $A_i$ words of weight $i$.
We know from Theorem~\ref{L3} that $E$ is closed under
componentwise multiplication.

Therefore the $A_2 + A_3$ words of weights 2 and 3 in $E$ must be disjoint,
or else $E$ would contain a word of weight 1.
Omitting these words from $E$, we are left with a code
of length $2^{m-1} - 2A_2 - 3A_3$,
dimension $\ge 2^{m-1} - m - A_2 - A_3$, and minimal distance 4.
This violates the optimality of shortened Hamming codes unless
$A_2 = A_3 =0$ and $E$ is itself
an extended
Hamming code of length $2^{m-1}$.
For $m \ge 5$ we now use \eqn{E30} to deduce that $E$
contains a word of weight 1, a contradiction. \hfill $\blacksquare$

Theorem~\ref{XH} demonstrates that a binary code can be $\ZZ_4$-linear, even though its dual is not.
For $RM(1,m)$ is $\ZZ_4$-linear, while in general its dual,
$RM(m-2,m)$, is not.
\subsection{The Golay code.~~}
Since the \NR code is $\ZZ_4$-linear
(as we shall see in Theorem~\ref{NR}) and is closely connected
with the Golay code (\cite[p.~73]{MS77},
it is natural to ask if the Golay code itself is $\ZZ_4$-linear.

\begin{theo}
\label{G}
The $[24,12,8]$ Golay code $G$ is not $\ZZ_4$-linear.
\end{theo}

\noindent{\em Proof}.~
Suppose on the contrary that $G$ is the binary image of a quaternary linear code $\sG$.
The swap map $\sigma$ (see \eqn{E1c}) is a fixed-point-free
involution in $\mbox{Aut} (G)$, the Mathieu group $M_{24}$.
It is known (\cite{ATLAS}, \cite{CS92}) that $M_{24}$ contains a single
conjugacy class of such involutions.
Therefore, without loss of generality, we may suppose that
this involution is the map defined by addition
of the hexacodeword $11\, \omega \, \omega \, \overline{\omega}\, \overline{\omega}$
in the MOG
description of $G$ (see \cite{CS92}, Chap.~11, \S9).
In the MOG diagram this is the permutation \\

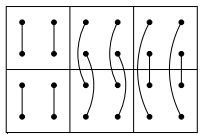
\vsp
\noindent
\centerline{\psfig{file=L2.ps,width=1in}}

\vsp
\noindent
The diagram specifies the division of the 24 coordinates into twelve pairs,
although we do not yet know which coordinate of each pair is on the left (in \eqn{E1b}) and which is on the right.
Consider the Golay codewords
$$
u = \begin{array}{|c|c|c|} \hline
11 & 11 & 00 \\
11 & 11 & 00 \\ \hline
00 & 00 & 00 \\
00 & 00 & 00 \\ \hline
\end{array}~,~~~~
v = \begin{array}{|c|c|c|} \hline
01 & 11 & 11 \\
10 & 00 & 00 \\ \hline
10 & 00 & 00 \\
10 & 00 & 00 \\ \hline
\end{array} ~.
$$
Then
$$
u + \sigma (u) = \begin{array}{|c|c|c|} \hline
00 & 11 & 00 \\
00 & 11 & 00 \\ \hline
00 & 11 & 00 \\
00 & 11 & 00 \\ \hline
\end{array} ~, ~~~~
v+ \sigma (v) = \begin{array}{|c|c|c|} \hline
11 & 11 & 11 \\
11 & 00 & 00 \\ \hline
00 & 11 & 00 \\
00 & 00 & 11 \\ \hline
\end{array}
$$
and $(u+ \sigma (u)) \ast (v+ \sigma (v))$ (which by Theorem~\ref{L3}
must be in $G$) has weight 4, a contradiction.~~\hfill $\blacksquare$
\section{Cyclic codes over $\ZZ_4$ and Galois rings}
\subsection{Galois rings.~~}
To study BCH and other cyclic codes of length $n$ over an alphabet of size $q$,
it is customary to work in a Galois field $GF(q^m)$, an extension of degree $m$ of a ground field $GF(q)$ \cite{MS77}.
The ground field $GF(q)$ is identified with the alphabet, and the extension field is chosen so that it contains an $n$th root of unity.

A similar approach is used for cyclic codes of length $n$ over $\ZZ_4$,
only now one constructs a Galois ring $GR (4^m)$ (not a field), that is an extension of $\ZZ_4$ of degree $m$ containing an $n$th root of unity.

Galois rings have been studied by MacDonald \cite{Mac74}, Liebler and Mena \cite{LM88},
Shankar \cite{Sha}, Sol\'{e} \cite{So89},
Yamada \cite{Ya90},
Bozta\c{s}, Hammons and Kumar \cite{BHK92}, among others, and of course the
general machinery of commutative algebra, as described for example in
Zariski and Samuel \cite{ZS60},
is applicable to these rings.
We list here some of the basic facts we shall need;
proofs may be found in the above references.

Let $h_2 (X) \in \ZZ_2 [X]$ be a primitive irreducible polynomial of degree $m$.
There is a unique monic polynomial $h(X) \in \ZZ_4 [X]$ of degree $m$
such that $h(X) \equiv h_2 (X)$ $( \bmod~2 )$ and $h(X)$
divides $X^n -1$ $( \bmod~4)$, where $n=2^m -1$
(see for example Yamada \cite{Ya90}).
The polynomial $h(X)$ is a {\em primitive basic irreducible} polynomial,
and may be found as follows.

Let $h_2 (X) = e(X) - d(X)$, where $e(X)$
contains only even powers and $d(X)$ only odd powers.
Then $h(X)$ is given by $h(X^2) = \pm (e^2 (X) - d^2 (X))$.
This is Graeffe's method \cite{Us48}, \cite{So89} for
finding a polynomial whose roots are the squares of the roots of
$h_2 (X)$.
For example, when $m=3$, $n=7$ we may take $h_2 = X^3 +X+1$.
Then $e=1$, $d= -X^3 -X$, $e^2 - d^2 = - X^6 - 2X^4 - X^2 +1$, so
\beql{E3.1}
h(X) = X^3 + 2X^2 + X-1 ~.
\eeq
Table~I in \cite{BHK92} gives all primitive basic irreducible
polynomials of degree $m \le 10$.

Let $\xi$ be a root of $h(X)$, so that $\xi^n =1$.
Then the {\em Galois ring} $GR (4^m)$ is defined to be $R = \ZZ_4 [ \xi ]$.
There are two canonical ways to represent the $4^m$ elements of $R$ (just as there are two canonical ways, multiplicative and additive, to represent elements
of $GF(q^m)$).

In the first representation, every element $c \in R$ has a unique
`multiplicative' or 2-adic representation
\beql{E3.2}
c = a+ 2b ~,
\eeq
where $a$ and $b$ belong to the set
\beql{E3.2a}
\sT = \{ 0,1, \xi , \xi^2 , \ldots , \xi^{n-1} \} ~.
\eeq
The map $\tau : c \mapsto a$ is given by
\beql{E3.3}
\tau (c) = c^{2^m} , ~~~c \in R ~,
\eeq
and satisfies
\begin{eqnarray}
\label{E3.4}
\tau (cd) & = & \tau (c) \tau (d) ~, \\
\label{E3.5}
\tau (c+d) & = & \tau (c) + \tau (d) + 2(cd)^{2^{m-1}}
\end{eqnarray}
(see \cite{Ya90}).
Given $c$, one determines $a$ from \eqn{E3.3} and
then $b$ from \eqn{E3.2}.

In the second representation, each element $c \in R$ has a unique
`additive' representation
\beql{E3.6}
c = \sum_{r=0}^{m-1} b_r \xi^r , ~~~
b_r \in \ZZ_4 ~.
\eeq
For example, if $m=3$ and $h$ is given by
\eqn{E3.1}, the additive representations for the elements of $\sT$ and $2 \sT$ are
\beql{E3.7}
\begin{array}{cc@{~}c@{~}cc@{~}c@{~}c}
\mbox{element} & b_0 & b_1 & b_2 & 2b_0 & 2b_1 & 2b_2 \\
~~ \\ [-.1in]
0 & 0 & 0 & 0 & 0 & 0 & 0 \\
1 & 1 & 0 & 0 & 2  & 0 & 0 \\
\xi & 0 & 1 & 0 & 0 & 2 & 0 \\
\xi^2 & 0 & 0 & 1 & 0 & 0 & 2 \\
\xi^3 & 1 & 3 & 2 & 2 & 2 & 0 \\
\xi^4 & 2 & 3 & 3 & 0 & 2 & 2 \\
\xi^5 & 3 & 3 & 1 & 2 & 2 & 2 \\
\xi^6 & 1 & 2 & 1 & 2 & 0 & 2
\end{array}
\eeq
This table may be produced (just as for Galois fields) by a
(modulo~4) shift
register whose feedback polynomial is $h(X)$.
By using \eqn{E3.2}, the table gives the additive
representation of every element of $R$.

One essential difference between $R= GR (4^m)$ and a Galois field is that $R$
contains zero divisors:
these are the elements of the radical $2R$, the unique maximal ideal in
$R$ ($R$ is a local ring).
Let
$\mu$ denote the map $R \to R/2R$.
Then $\theta = \mu ( \xi )$ is
a root of $h_2 (X)$, and we can identify $R/2R$ with $GF(2^m)$, taking the elements of $GF(2^m)$ to be
\beql{E3.8}
\mu ( \sT) = \{ 0, 1, \theta , \theta^2 , \ldots , \theta^{n-1} \} ~.
\eeq
We denote the set of {\em regular} or {\em invertible} elements of $R$ by
$R^\ast = R~\setminus~2R$.
Every element of $R^\ast$ has a unique representation in the form
$\xi^r (1+2t)$, $0 \le r \le n-1$, $t \in \sT$.
$R^\ast$ is a multiplicative group of order $(2^m -1)2^m$
which is a direct product $H \times \sE$, where $H$ is the cyclic
group of order $2^m -1$ generated by $\xi$, and $\sE$
is the group of {\em principal units} of $R$, that is,
elements of the form $1+ 2t$, $t \in \sT$.
$\sE$ has the structure of an elementary abelian group of order $2^m$
and is isomorphic to the additive group of $GF(2^m)$.

\subsection{Frobenius and trace maps.~~}
The {\em Frobenius map} $f$ from $R$ to $R$ is the ring automorphism
that takes any element
$c = a+ 2b \in R$ to
\beql{E3.9}
c^f = a^2 + 2b^2 ~.
\eeq
$f$ generates the Galois group of $R$ over $\ZZ_4$, and $f^m$ is the identity map.
The {\em relative trace} from $R$ to $\ZZ_4$ is defined by
\beql{E3.10}
T(c) = c + c^f + c^{f^2} + \cdots + c^{f^{m-1}}, ~~~c \in R ~.
\eeq
For comparison, the usual trace from $GF(2^m)$ to $\ZZ_2$ is given by
\beql{E3.11}
tr (c) = c+c^2 + c^{2^2} + \cdots + c^{2^{m-1}} , ~~~
c \in GF(2^m) ~,
\eeq
and the \For map is simply the squaring map
\beql{E3.11a}
f_2 (c) = c^2 , ~~~c \in GF(2^m) ~.
\eeq
The following commutativity relationships between
these maps
are easily verified:
\bary
\label{E3.11b}
\mu \circ f & = & f_2 \circ \mu ~, \\
\label{E3.11c}
\mu \circ T & = & tr \circ \mu ~.
\eary
In particular, since $tr$ is not identically zero, it follows that the Galois
ring trace is nontrivial.
In fact, $T$ is an onto mapping from $R$ to $\ZZ_4$.
The set of elements of $R$ invariant under $f$ is identical with $\ZZ_4$.
\subsection{Dependencies among $\xi^j$.~}
For later use we record some results about dependencies among the powers
$\xi^j$.

(P1)~$\pm \xi^j \pm \xi^k$ is invertible for $0 \le j < k < 2^m -1$,
for $m \ge 2$.
Proof.~If on the contrary we had $\pm \xi^j \pm \xi^k = 2 \la$, $\la \in R$,
then applying $\mu$ we obtain $\theta^j + \theta^k =0$,
which contradicts the fact that $\theta$ is primitive in $GF(2^m)$. \hfill $\blacksquare$

\vsp
(P2)~$\xi^j - \xi^k \neq \pm \xi^l$ for
distinct $j,k,l$ in the range
$[0, 2^m -2]$, for $m \ge 2$.
Proof.~Otherwise, after rearranging, we have $1+ \xi^a = \xi^b$ for $a \neq b$.
Squaring gives $1+2 \xi^a + \xi^{2a} = \xi^{2b}$, but applying
the \For map gives $1+ \xi^{2a} = \xi^{2b}$, so
$2 \xi^a = 0$, a contradiction. \hfill $\blacksquare$

\vsp
(P3)~Suppose $i,j,k,l$ are in the range
$[0,2^m -2]$ and
$i \neq j$,
$k \neq l$, $m \ge 3$.
Then
$$\xi^i - \xi^j = \xi^k - \xi^l
\Leftrightarrow i=k ~~~\mbox{and}~~~j=l ~,
$$
Proof.~Suppose $1+ \xi^a = \xi^b + \xi^c$.
Squaring and subtracting the result of applying the \For
map gives $2 \xi^a = 2 \xi^{b+c}$.
Therefore $\xi^a \equiv \xi^{b+c}$ ($\bmod~2$), so if we write
$x= \theta^a$, $y= \theta^b$, $z= \theta^c$ we have
$x = yz$.
But also $1+ x = y+z$, so $(y+1) (z+1) =0$, which since $\theta$ is primitive
in $GF(2^m)$ implies $y$ or $z=1$.  \hfill $\blacksquare$

\vsp
(P4)~For odd $m \ge 3$,
$$\xi^i + \xi^j + \xi^k + \xi^l = 0 \Rightarrow
i=j=k=l ~.$$
Proof.~Suppose $\xi^a + \xi^b + \xi^c =-1$.
Arguing as in the previous proof we obtain
$x^2 + y^2 = (x+1) (y+1)$,
hence $u^2 + v^2 =uv$,
with $x = u+1$, $y= v+1$.
Substituting $v= tu$ we find
$u^2 (t^2 + t+1)=0$.
But $t^2 + t+1 \neq 0$ in $GF(2^m)$, $m$ odd,
since $tr (t^2 + t+1) = m \neq 0$,
so $u=0$, $x=1$, therefore $a= b=c=0$.  \hfill $\blacksquare$

\vsp
Properties P2, P3 and P4 are also consequences of the fact that errors
of weight $\le 2$ in the `Preparata' code can be decoded uniquely, as shown in
\S5.3.
\subsection{The ring $\sR$.~~}
As usual when studying cyclic codes of length $n$ it is convenient
to represent codewords by polynomials modulo $X^n -1$.
We identify $v= (v_0 , v_1 , \ldots , v_{n-1} )$ with the polynomial
$v(X) = \sum_{r=0}^{n-1} v_r X^r$ in the ring
$\sR = \ZZ_4 [X] / (X^n -1)$.
We must be careful when working with $\sR$:
it is not a unique factorization domain --- for example $X^4 -1$ has two
distinct factorizations into irreducible polynomials in $\sR$:
\bara
X^4 -1 & = &
(X-1) (X+1)(X^2+1) \\
& = & (X+1)^2 (X^2 + 2X-1) ~.
\eara
Note also that every element $1+ 2 \la$, $\la \in R$,
is a root of $X^2 -1$.
On the other hand $\sR$ {\em is} a principal ideal domain:
just as in the binary case, cyclic codes have a single generator
(the proof is given in Calderbank and Sloane,
Modular and $p$-adic cyclic codes,
{\em Designs, Codes and Cryptography}, to appear).
\section{Kerdock codes}
\hsp
The main result of this section is a very simple quaternary
construction for \Ker codes.
\subsection{The Kerdock code is an extended cyclic code over $\ZZ_4$.~~}
Let $h(X)$ be a primitive basic irreducible polynomial of degree $m$,
as above, and let $g(X)$ be the reciprocal polynomial to
$(X^n -1 ) / ((X-1) h(X))$,
where $n=2^m -1$.
\begin{theo}
\label{K}
let $\sK^-$ be the cyclic code of length $n$ over $\ZZ_4$
with generator polynomial $g(X)$, and let $\sK$
be obtained from $\sK^-$ by adjoining a zero-sum
check symbol.
Then for odd $m \ge 3$ the binary image $K = \phi ( \sK )$ of $\sK$ under the
Gray map \eqn{E1} is a nonlinear code of length $2^{m+1}$, with
$4^{m+1}$ words and minimal distance $2^m -2^{(m-1)/2}$ that is
equivalent to the \Ker code.
This code is distance invariant.
\end{theo}

Note that $\sK^-$ has parity check polynomial $(X-1) h(X)$.
There are two equivalent generator matrices for $\sK$.
The first is
\beql{E3.12}
\left[
\begin{array}{cccccc}
1 & 1 & 1 & 1 & \cdots & 1 \\
0 & 1 & \xi & \xi^2 & \cdots & \xi^{n-1}
\end{array}
\right] ~,
\eeq
where the entries in the second row are to be replaced by the
corresponding $m$-tuples $(b_0 b_1 \cdots b_{m-1} )'$ (the prime
indicating transposition) obtained from \eqn{E3.6}.
Alternatively,
let $g(X) = \sum_{j=0}^\delta g_j X^j$,
$\delta =2^m -m-2$,
$g_j \in \ZZ_4$,
and let $g_\infty = - \sum_{j=0}^\delta g_j$.
Then the second form for the generator matrix for $\sK$ is
\beql{E3.13}
\left[
\begin{array}{cccccccc}
g_\infty & g_0 & g_1 & \cdots & g_\delta & 0 & \cdots & 0 \\
g_\infty & 0 & g_0 & \cdots & g_{\delta-1} & g_\delta & \cdots & 0 \\
\cdot & \cdot & \cdot & ~ & \cdot & \cdot & ~ & ~ \\
g_\infty & 0 & 0 & \cdots & g_0 & g_1 & \cdots & g_\delta
\end{array}
\right] ~.
\eeq
$\sK$ is a code of type $4^{m+1}$.
The binary code $K^{(1)}$ associated with $\sK$ (see \eqn{E33}) is $RM(1,m)$.

For example, with $m=3$ and $h$ given by \eqn{E3.1},
we find $g= x^3 + 2x^2 +x-1$, so the two equivalent
generator matrices are
\beql{E3.14}
\left[
\begin{array}{cccccccc}
1 & 3 & 1 & 2 & 1 & 0 & 0 & 0 \\
1 & 0 & 3 & 1 & 2 & 1 & 0 & 0 \\
1 & 0 & 0 & 3 & 1 & 2 & 1 & 0 \\
1 & 0 & 0 & 0 & 3 & 1 & 2 & 1
\end{array}
\right] ~, ~~~
\left[
\begin{array}{cccccccc}
1 & 1 & 1 & 1 & 1 & 1 & 1 & 1 \\
0 & 1 & 0 & 0 & 1 & 2 & 3 & 1 \\
0 & 0 & 1 & 0 & 3 & 3 & 3 & 2 \\
0 & 0 & 0 & 1 & 2 & 3 & 1 & 1
\end{array}
\right]
\eeq
(the second one being read from \eqn{E3.7}).

For $m=5$,
we may take $h(X) = \sum_{i=0}^5 h_i X^i$,
$g(X) = \sum_{i=0}^{25} g_i X^i$,
where $h_0 \ldots$ and $g_0 \ldots$ are
323001 and 11120122010303133013212213.

\Ker codes contain more codewords than any known linear code with the same minimal distance (although we are not aware of any theorem to guarantee
this, except at length 16).
\subsection{Family $\sA$.~~}
If we omit the factor $X-1$ from the parity check polynomial
for $\sK^-$, we obtain a cyclic code containing $4^m$
codewords.
Let $\sA$ denote the family of cyclically distinct vectors obtained from this code by
deleting the zero vector and failing to distinguish between a vector and
any of its cyclic shifts.
The corresponding collection $\Omega ( \sA )$ of complex-valued sequences
has been studied in \cite{Bo90}, \cite{BHK92}, \cite{So89}, \cite{US91} as a family
of asymptotically optimal CDMA signature sequences
(referred to as Family $\sA$ in \cite{BHK92}).
Since the sequences of $\Omega ( \sA )$ have low values of auto- and cross-correlation, the set $\Omega ( \sA )$ also has large minimal Euclidean distance.
\subsection{Trace description of Kerdock code and proof of Theorem~\protect\ref{K}.~~}
\begin{theo}
\label{KT}
The codes $\sK^-$ and $\sK$ have the following trace
descriptions
over the ring $R$.
\begin{itemize}
\item[(a)]
$c = (c_0 , c_1 , \ldots , c_{n-1} )$ is a codeword in $\sK^-$ if
and only if, for some $\la \in R$ and $\epsilon \in \ZZ_4$,
\beql{E.KT.1}
c_t = T( \la \xi^t ) + \epsilon , ~~~
t \in \{ 0,1, \ldots , n-1 \} ~.
\eeq
Thus
\beql{E3.15}
\sK^- = \{ \epsilon~ \bone + v^{( \la )} : ~
\epsilon \in \ZZ_4 , ~~~\la \in R \} ~,
\eeq
where
$$
v^{( \la )} = (T( \la ), T( \la \xi ) , T( \la \xi^2) , \ldots , T( \la \xi^{n-1} )) ~.
$$
\item[(b)]
$c = (c_\infty , c_0 , c_1 , \ldots , c_{n-1} )$ is a codeword in
$\sK$ if and only if, for some $\la \in R$ and $\epsilon \in \ZZ_4$,
\beql{E.KT.2}
c_t = T( \la \xi^t ) + \epsilon , ~~~t \in \{ \infty , 0,1, \ldots , n-1 \}
~,
\eeq
with the convention that $\xi^\infty =0$.
\end{itemize}
\end{theo}

This theorem is essentially equivalent to Theorem~3 of
\cite{BHK92}.

\PF
(a)~Let $C$ be the code defined by \eqn{E3.15}.
If $c(X)$ is
the polynomial form of a
codeword in $C$, then $c(X) (X-1) h(X) =0$ [the all $1$'s
vector is annihilated by $X-1$ and the $v^{( \la )}$ by $h(X)$].
Therefore $C \subseteq \sK^-$.
Since $C$ and $\sK^-$
contain the same number of codewords, $C= \sK^-$.
(b)~follows because the zero-sum check for $\epsilon \, \bone$ is
$\epsilon$ and for $v^{( \la )}$ it is 0. \hfill $\blacksquare$

\vsp
\noindent{\em Proof of Theorem~\ref{K}}.
We consider an arbitrary codeword $c \in \sK$ in the form
\eqn{E.KT.2}.
We will show that $c_t$ has 2-adic expansion
\beql{E3.16a}
c_t = a_t + 2b_t , ~~~
t \in \{ \infty , 0,1, \ldots , n-1 \} ~,
\eeq
given by
\bary
\label{E3.16}
a_t & = & tr ( \pi \theta^t ) + A ~, \\
b_t & = & tr ( \eta \theta^t ) + Q ( \pi \theta^t ) + B ~, \label{E3.17}
\eary
where the elements $\pi$, $\eta \in GF(2^m)$ and $A,B \in \ZZ_2$ are
arbitrary,
$$Q(x) = \sum_{j=1}^{(m-1)/2} tr (x^{1+2^j} ) ~,~~~
x \in GF (2^m) ~,
$$
and we adopt the convention that $\theta^\infty = 0$.

Let $\la = \xi^r + 2 \xi^s$,
$r$, $s \in \{ \infty , 0 , \ldots , n-1 \}$, so that
$$c_t = \epsilon + T ( \xi^{r+t} ) + 2T ( \xi^{s+t} ) =
a_t + 2b_t ~.$$
Projecting modulo 2, we obtain
$$a_t = \af ( \epsilon ) + tr ( \pi \theta^t ) ~,$$
where $\pi = \mu ( \xi^r )$, $\theta = \mu ( \xi )$.
To find $b_t$, we compute
$c_t - c_t^2 = 2b_t$
(since $a_t = 0$ or 1) and obtain
\begin{eqnarray*}
2b_t & = & ( \epsilon - \epsilon^2) +
(T( \xi^{r+t}) - T^2 (\xi^{r+t}))+ 2 \epsilon T ( \xi^{r+t})
+ 2T( \xi^{s+t} ) \\
& = & 2 \beta ( \epsilon ) + 2 \sum_{0 \le j < k \le m-1}
( \xi^{r+t})^{2^j + 2^k} +
2T (( \epsilon \xi^r + \xi^s ) \xi^t ) ~.
\end{eqnarray*}
Thus
$$
b_t = \beta ( \epsilon ) + Q ( \pi \theta^t) + tr ( \eta \theta^t ) ~,$$
where $\eta = \mu ( \epsilon \xi^r + \xi^s )$.

The next step
is to observe that the vectors $(b_t)$ and $(a_t + b_t)$ defined by
\eqn{E3.16}, \eqn{E3.17} are the left and right halves of the codewords
in Kerdock's original definition
(\cite{Ke72}; \cite[p.~458]{MS77}).
But the Gray map $\phi$ sends $c$ to $( \beta (c) , \gamma (c)) = ((b_t)$,
$(a_t + b_t))$.

The fact that $\phi ( \sK )$ is distance invariant follows from
Theorem~\ref{tH2.1}. \hfill $\blacksquare$

It is shown in \cite{BK93} that when $m$ is odd, the family of binary
sequences
$\{Q( \pi \theta^t ) + tr ( \eta \theta^t ): \eta, \pi$ in
$GF (2^m )$, not both zero$\}$ has
Gold-like correlation properties, but a larger linear span.
\subsection{The first-order Reed-Muller subcode.~~}
The vectors for which $\pi =0$ in
\eqn{E3.16}, \eqn{E3.17} form a linear subcode of $\sK$, with
generator matrix
$$
\left[
\matrix{1 & 1 & 1 & 1 & \cdots & 1 \cr
0 & 2 & 2 \xi & 2 \xi^2 & \cdots & 2 \xi^{n-1} \cr}
\right] ~,
$$
whose binary
image is the first-order Reed-Muller code contained in the \Ker code.
\subsection{The Nordstrom-Robinson code.~~}
The case $m=3$ is particularly interesting.
The Kerdock and Preparata codes of length 16 coincide, giving the \NR
code (\cite{NR67}; see also \cite{SeZ}).
This is the unique binary code of length 16, minimal distance 6,
containing 256 words \cite{Sn73}, \cite{Go77}.
In this case $\sK$ is the `octacode', whose generator matrix is given in
\eqn{E3.14}.
The octacode may also be characterized as the unique
self-dual quaternary code of length 8 and minimal Lee weight 6 \cite{CS93},
or as the `glue code' required to construct the 24-dimensional Leech lattice from eight copies
of the face-centered cubic lattice \cite[Chap.~24]{CS92}.
Thus the following theorem is a special case of Theorem~\ref{K}.

\begin{theo}
\label{NR}
The \NR code is the binary image of the octacode under the Gray map.
\end{theo}

\noindent
The symmetrized weight enumerator of the octacode is (\cite{CS93})
$$W^8 + 16X^8 + Y^8 + 14W^4 Y^4 + 112WX^4 Y(W^2 + Y^2) ~,$$
and the weight distribution of the \NR code is then given by \eqn{E2.14}.
\subsection{Weight distribution.~~}
The weight distribution of any \Ker code is also easily determined
from the new quaternary description.
\begin{theo}
\label{KW1}
The binary \Ker code $K = \phi ( \sK )$ of length $2^{m+1}$ ($m$ odd $\ge 3$)
has the following weight distribution:
\beql{E3.49}
\begin{array}{cc}
i & A_i \\
~~ \\ [-.1in]
0 & 1 \\
~~ \\ [-.05in]
2^m - 2^{(m-1)/2} & 2^{m+1} (2^m -1) \\
~~ \\ [-.05in]
2^m & 2^{m+2} -2 \\
~~ \\ [-.05in]
2^m + 2^{(m-1)/2} & 2^{m+1} (2^m -1) \\
~~ \\ [-.05in]
2^{m+1} & 1
\end{array}
\eeq
(cf. \cite{MS77}, Fig.~15.7).
\end{theo}
\PF
This is a slight modification of the argument used in \cite{BHK92} to obtain the correlation distribution of the associated complex sequences.
We assume the codewords $c \in \sK$ are defined
as in Theorem~\ref{KT}.
As mentioned in \S4.5,
the words for which $\pi =0$ (and $\la \not\in R^\ast$ in \eqn{E3.15})
form a first-order Reed-Muller code, and account for the words of
weights 0, $2^m$ and $2^{m+1}$.

We now consider a word $v^{( \la )} \in \sK^-$ for $\la \in R^\ast$.
Let $n_j = n_j (v^{( \la )})$ (see \eqn{E2.3}).
We claim that there exist $\delta_1$, $\delta_2 = \pm 1$ so that
\beql{E3.50}
\begin{array}{ll}
n_0 = 2^{m-2} -1 + \delta_1 2^{(m-3)/2} , & n_1 = 2^{m-2} + \delta_2 2^{(m-3)/2} ~, \\
~~ \\
n_2 = 2^{m-2} - \delta_1 2^{(m-3)/2} , & n_3 = 2^{m-2} - \delta_2 2^{(m-3)/2} ~.
\end{array}
\eeq
Let
$$S = \sum_{j=0}^{2^m -2} i^{T( \la \xi^j)} =
n_0 - n_2 + i ( n_1 - n_3 ) ~.
$$
Then
$$|S|^2 =
2^m -1 + \sum_{j \neq k} i^{T( \la ( \xi^j - \xi^k))} ~.$$
We use properties (P1), (P2), (P3) to rewrite this as
$$|S|^2 = 2^m -1 + \sum_{\nu \in R^\ast} i^{T( \nu )} - S - \overline{S} ~.$$
But it is easily verified that
$$
\sum_{\nu \in R^\ast} i^{T( \nu )} = 0$$
(see \cite[p.~1104]{BHK92}), hence
$$
\begin{array}{c}
(S+1) ( \overline{S} +1) =2^m ~, \\
~~ \\
(n_0 - n_2 + 1 )^2 + (n_1 - n_3)^2 = 2^m ~.
\end{array}
$$
The diophantine equation $X^2 + Y^2 =2^m$ has a
unique solution, so
\bary
\label{E3.51}
n_0 - n_2 & = & -1 \pm2^{(m-1)/2}~, \\
\label{E3.52}
n_1 - n_3 & = & \pm 2^{(m-1)/2} ~.
\eary
We also know that $\mu ( v^{( \la )})$ is in the simplex code, so
\bary
\label{E3.53}
n_1 + n_3 & = & 2^{m-1} ~, \\
\label{E3.54} n_0 + n_2 & = & 2^{m-1} -1 ~.
\eary
\eqn{E3.50} follows from \eqn{E3.51}--\eqn{E3.54}.

We now consider the four words of $\sK$ obtained from $\epsilon \bone + v^{( \la )}$
$( \epsilon =0,1,2,3 )$ by appending the zero-sum check symbol $\epsilon$.
For $\bone + \vla$, for example, we have
$$
\begin{array}{ll}
n_1 = 2^{m-2} + \delta_1 2^{(m-3)/2} , & n_2 = 2^{m-2} + \delta_2 2^{(m-3)/2} ~, \\
~~ \\
n_3 = 2^{m-2} - \delta_1 2^{(m-3)/2} , & n_0 = 2^{m-2} - \delta_2 2^{(m-3)/2} ~,
\end{array}
$$
which is
a word of Lee weight
$$n_1 + n_3 + 2n_2 = 2^m + \delta_2 2^{(m-1)/2} ~.$$
Of these four words obtained from $\vla$, two have
Lee weight $2^m + 2^{(m-1)/2}$ and two have Lee weight $2^m - 2^{(m-1)/2}$.
This holds for all $2^m (2^m -1 )$ words $\vla$,
$\la \in R^\ast$, and establishes \eqn{E3.49}.~~\hfill $\blacksquare$

\vsp
When $m$ is even, $m \ge 2$, a similar argument shows that $\phi ( \sK )$ is a nonlinear code of length $2^{m+1}$, with $4^{m+1}$ codewords,
minimal distance $2^m - 2^{m/2}$, and weight distribution
$$
\begin{array}{cc}
i & A_i \\
~~ \\ [-.1in]
0 & 1 \\
~~ \\ [-.05in]
2^m - 2^{m/2} & 2^m ( 2^m -1 ) \\
~~ \\ [-.05in]
2^m & 2^{m+1} (2^m +1 ) -2 \\
~~ \\ [-.05in]
2^m + 2^{m/2} & 2^m (2^m -1) \\
~~ \\ [-.05in]
2^{m+1} & 1
\end{array}
$$
This code is not as good as a double-error-correcting BCH code.
\subsection{Soft-decision decoding of Kerdock codes.~~}
Although in the theoretical development we make a distinction
between the quaternary code $\sK$ and the associated nonlinear binary code
$K = \phi ( \sK )$ (and similarly in Section~V between $\sP = \sK^\perp$
and $P = \phi ( \sP ))$, they are really two different descriptions of the same code.
For instance, a decoder for the quaternary code obviously provides
a decoder for the binary code and conversely.

The following is a new soft-decision decoding algorithm for the \Ker code.
This is comparable in complexity to previously known techniques
that were
derived from the binary description of the code.

The idea is to extend the fast Hadamard transform (FHT)
soft-decision decoding algorithm for the binary first-order Reed-Muller
code to the Kerdock code.
This provides substantial savings over brute-force correlation decoding.
Define
$$\Delta = \{ \infty , 0,1,2, \ldots , n-1 \} ,~~n=2^m-1~.$$
Brute-force decoding of a received vector
$\{ v_t: t \in \Delta \}$ requires the computation of its correlation with all possible received signals.
In particular, the decoder must compute the correlation
$$\zeta ( \la , \delta ) = \sum_{t \in \Delta} v_t i^{-[T ( \la \xi^t ) + \delta ]}$$
for all $\la = \xi^r + 2 \xi ^s$, $r,s \in \Delta$ and all
$\delta \in \ZZ_4$, and find that pair $( \la , \delta )$ for which
$\mbox{Real} \{ \zeta ( \la , \delta ) \}$ is a maximum.
Computed directly, this technique requires
$4^{m+1} 2^m$ multiplications and
\linebreak
$4^{m+1} (2^m -1)$ additions.

An immediate reduction in complexity is obtained by writing
$$\zeta ( \la , \delta ) = i^{- \delta}
\sum_{t \in \Delta} v_t i^{-T( \xi^{t+r})} (-1)^{tr (\theta^{t+s})}~,$$
where we adopt the convention that for $l$ in $\Delta$,
$l+ \infty = \infty$.
The correlation sums
$\zeta ( \xi^r + 2\xi^s , \delta )$ may now be viewed
(after some reordering of indices) as $i^{- \delta}$ times the
Hadamard transform of the $2^m$ complex vectors
$\{ v_t i^{-T( \xi^{t+r})} \}$ of length $2^m$.
Using the FHT, each of these can be computed using
$m 2^m$ additions/subtractions.
Thus the overall requirement is for about $4^m$
multiplications (one multiplicand is always a power of $i$)
and $m 4^m$ additions/subtractions.

This complex-data FHT decoding algorithm is of the same order of complexity as recently published real-data FHT decoders for the \Ker codes
\cite{Ad87}, \cite{ELN88}
based on the general super-code decoding method of Conway and Sloane \cite{CS86}.
These real-data algorithms perform $2^m$ FHTs of size $2^{m+1}$.
Finally, we note that the case $m=3$ corresponds to decoding the Nordstrom-Robinson code.
\section{Preparata codes}
\hsp
In this section we show that the binary image of the dual code $\sP = \sK^\perp$ is a Preparata-like code with
essentially the same properties as Preparata's original code (yet
is much simpler to construct).
\subsection{The `Preparata' code is an extended cyclic code over $\ZZ_4$.~~}
Let $h(X)$ and $g(X)$ be defined as in \S4.1.

\begin{theo}
\label{P}
Let $\sP^-$ be the cyclic code of length $n=2^m -1$ with generator
polynomial $h(X)$, and let $\sP$ be obtained from $\sP^-$ by adjoining
a zero-sum check symbol, so that $\sP = \sK^\perp$.
Then for odd $m \ge 3$ the binary image $P = \phi ( \sP )$ of $\sP$ under the
Gray map \eqn{E1} is a nonlinear code of length $l=2^{m+1}$,
with $2^{l-2m-2}$ codewords and minimal distance 6.
This code is distance invariant and its weight distribution
is the MacWilliams transform of the weight
distribution of the \Ker code of the same length.
\end{theo}

Note that $\sP^-$ has parity check polynomial $g(X)$,
and that \eqn{E3.12}, \eqn{E3.13} are equivalent parity check matrices for $\sP$.
Also $\sP$ is a code of type $4^{2^m -m-1}$.
The code $P = \phi ( \sP )$ is the $\ZZ_4$-dual of $K$, and we refer
to it as a `Preparata' code, using the quotes to distinguish it from
Preparata's original code.
It is known that the \Pre code (and $P$) contains more codewords than any
linear code with the same minimal distance \cite{BT92}.
The binary code $P^{(1)}$ associated with $\sP$ (see \eqn{E33}) is
$RM (m-2, m)$.

\noindent{\em Proof of Theorem~\ref{P}}.
It follows from Theorems~\ref{tH2.1} and \ref{H2.5} that $P$ is distance invariant and its weight distribution is the MacWilliams transform
of that of $K$.
By Theorem~24 of \cite{MS77}, Chapter~15, $P$ has the same weight distribution as the original Preparata code. \hfill $\blacksquare$

Semakov, Zinoviev and Zaitsev \cite{SZZ} had already
shown in 1971 that any code with the same
parameters as the Preparata code must be distance invariant.

The decoding algorithm given below provides an alternative proof that $\sP$ has
minimal Lee weight 6, for odd $m \ge 3$.
For
even $m \ge 2$, $\sP$ contains words of Lee weight 4.
For $\xi$ satisfies $\xi^{3t} =1$,
where $t = (2^m -1) /3$, and since $\xi^t-1 \in R^\ast$,
by (P1), $\xi^{2t} + \xi^t + 1 =0$, yielding a word of Lee weight 3 in $\sP^-$.

There is one essential difference between $P$ and the original \Pre code.
It is known that the latter is contained in the extended Hamming
code spanned by its codewords.
\begin{theo}\label{PD}
For odd $m \ge 5$, $P$ is contained in a nonlinear code with the same
weight distribution as the extended Hamming code of the same length,
and the linear code spanned by the codewords of $P$ has
minimal weight $2$.
\end{theo}

\PF
The first assertion follows by considering the binary images of the following sequence of codes:
\beql{E111.1}
ZRM (1,m) \subseteq \sK
\subseteq ZRM (2,m) \subseteq \cdots \subseteq
ZRM (2,m)^\perp \subseteq \sP \subseteq
ZRM (1,m)^\perp ~.
\eeq

For the second assertion we use the fact that $\sP$ is an extended cyclic code with generator polynomial
$h(X) = \sum_{j=0}^m h_j X^j$ (say).
Let $h_\infty = - h (1) = \pm 1$, since $h_2 (1) =1$.
Then $\sP$ has a generator matrix of the form
\beql{E4.1}
\begin{array}{c}
a \\ ~ \\ ~ \\ b \\ ~
\end{array}
\left[
\begin{array}{cccccccc}
h_\infty & h_0 & h_1 & \cdots & h_m & ~ & ~ & ~ \\
h_\infty & 0 & h_0 & \cdots & ~ & h_m & ~ & ~ \\
\cdot & \cdot & \cdot & \cdots & ~ & \cdot & ~ & ~ \\
h_\infty & 0 & 0 & \cdots & 0 & h_0 & \cdots & h_m \\
\cdot & \cdot & \cdot & \cdots
\end{array}
\right] ~.
\eeq
It follows from \eqn{E3} that the linear span of $P = \phi ( \sP )$ contains all words of the form
$\phi (2 \af (a) \ast \af (b))$, for $a,b \in \sP$.
Taking $a$ and $b$ to be as indicated in \eqn{E4.1} produces a word of
weight 2 in the linear span. \hfill $\blacksquare$

Again there is a result of Zaitsev, Zinoviev and Semakov that is relevant:
they showed in \cite{ZZS} that any code with the same parameters as the Preparata code is a subcode of a possibly nonlinear code with the same parameter
as an extended Hamming code.
Theorem~\ref{PD} answers a
question raised in that paper, by providing an example where the
Hamming-type code is indeed nonlinear.

As a quaternary linear code, $ZRM (1,m)^\perp$ (see \eqn{E111.1}) is the union of
$2^m -1$ translates of $\sP$,
each nonzero translate having minimal Lee weight 4.
The codewords of weight 4 in the binary image of $ZRM (1,m)^\perp$ (a nonlinear code with the same parameters an extended Hamming code) form a
Steiner system $S(3,4,2^{m+1})$.
It is not difficult to show that this
$S(3,4,2^{m+1} )$ is identical to the Steiner system formed by the codewords
of weight 4 in the classical extended Hamming code of length
$2^{m+1}$.
The blocks of this design are divided equally among the binary images of the
$2^m -1$ nonzero cosets of $\sP$.
The blocks falling in the binary image of a fixed coset
form a Steiner system $S(2,4,2^{m+1} )$.
\subsection{Transform-domain characterization of `Preparata' codes.~~}
In spite of the previous theorem, in this section we shall show
that
the `Preparata' code $P = \phi ( \sP )$ and Preparata's original code
have similar characterizations by finite field transforms.

We define the {\em Galois ring transform}
$\widehat{c} = ( \widehat{c} ( \la ))$, $\la =0,1, \ldots , n-1 , n=2^m -1$,
of a quaternary sequence
$c= (c_t)$, $t=0,1, \ldots, n-1$, by
$$\widehat{c} ( \la ) =
\sum_{t=0}^{2^m-2} c_t \xi^{\la t} ~.$$
The inversion formula
$$c_t = - \sum_{\la =0}^{n-1}
\widehat{c} ( \la ) \xi^{- \la t}$$
follows in the usual way from the fact that
$$
\sum_{\la =0}^{n-1} \xi^\la =
\frac{1- \xi^n}{1- \xi}
= 0 ~.
$$

We define the {\em finite field transform}
$\tilde{a} = ( \tilde{a} ( \la ))$, $\la =0,1, \ldots , n-1$, of a
{\em binary} sequence
$a= (a_t)$, $t=0,1, \ldots , n-1$, by
$$\tilde{a} ( \la ) = \sum_{t=0}^{n-1} a_t \theta^{\la t} ~,$$
where $\theta \in GF(2^m)$ is the image of $\xi \in R$ after reduction modulo 2
(as in \S3.1).
We define the {\em half-convolution} $\sH ( \tilde{a} , \la ) \in GF (2^m)$
of the sequence
$\tilde{a}$ at $\mbox{lag}~ \la$ by
$$\sH ( \tilde{a} , \la ) = \sum_{\la_1 \le \la_2 \atop \la_1 + \la_2 = \la}
\tilde{a} ( \la_1) \tilde{a} ( \la_2) ~,
$$
where $\la_1$, $\la_2 =0,1, \ldots , n-1$.
The summation is a half rather than full convolution because we exclude the cases
$\la_1 > \la_2$.

\begin{theo}
\label{C16A}
The quaternary `Preparata' code $\sP$ consists of all vectors $c = (c_t ) \in \ZZ_4^n$,
$t \in \{ \infty , 0,1, \ldots , n-1 \}$ satisfying the Galois ring transform
constraints
\begin{eqnarray}
\label{E.A0}
c_\infty + \widehat{c} (0) & = & 0 ~, \nonumber \\
\widehat{c} (1) & = & 0 ~.
\end{eqnarray}
\end{theo}

\PF
This follows from the definition of the Galois ring transform and
the parity check matrix for $\sP$ given in Eq.~\eqn{E3.12}. \hfill $\blacksquare$

\begin{theo}
\label{GRT}
The binary `Preparata' code $P$ consists of all vectors
$(b, a+b)$ for which $a,b \in \ZZ_2^{n}$ satisfy
\bary
\label{E115.0}
\tilde{a} (0) + a_\infty & = & 0 ~, \nonumber \\
\tilde{a} (1) & = & 0 ~, \nonumber \\
\tilde{b} (0) + b_\infty & = & \sH ( \tilde{a},0) + a_\infty ~, \nonumber \\
\tilde{b} (1) & = & \sH ( \tilde{a} ,1) ~.
\eary
\end{theo}

Note that equations \eqn{E.A0} are over $R$, whereas equations
\eqn{E115.0} are
over $GF (2^m)$.

\PF
Consider a codeword $c = (c_t) \in \sP$, where
$c_t = a_t + 2b_t$, $t \in \{ \infty , 0 , \ldots , n-1 \}$.
It follows from the previous theorem that
\begin{eqnarray}
\label{E.A1}
a_\infty + 2b_\infty + \widehat{a} (0) + 2 \widehat{b} (0) & = & 0 ~, \nonumber \\
\widehat{a} (1) + 2 \widehat{b} (1) & = & 0 ~.
\end{eqnarray}
The next step is identify the constraints that \eqn{E.A1} places on
$\tilde{a} ( \la )$, $\tilde{b} ( \la )$.
Given
\linebreak
$\la \in \{ 0,1, \ldots , n-1 \}$, let
\beql{E.A2}
\widehat{a} ( \la ) = e_\la + 2 f_\la , ~~~\mbox{where}~~~
e_\la , ~f_\la \in \sT ~.
\eeq
We find $f_\la$ indirectly, starting from the inversion formula
$$a_t = - \sum_{\la =0}^{n-1} \widehat{a} ( \la ) \xi^{- \la t} ~.$$
After squaring and also applying the Frobenius map we obtain
\begin{eqnarray*}
a_t & = & a_t^2 ~=~ \sum_{\la =0}^{n-1}
e_\la^2 \xi^{-2 \la t}
+ 2 \sum_{\la_1 < \la_2} e_{\la_1} e_{\la_2}
\xi^{- ( \la_1 + \la_2) t} \\
\noalign{and}
a_t & = & - \sum_{\la =0}^{n-1}
(e_\la^2 + 2 f_\la^2 ) \xi^{- 2 \la t} \\
& = & \sum_{\la =0}^{n-1} e_\la^2 \xi^{- 2 \la t} +
2 \sum_{\la =0}^{n-1} ( e_\la^2 + f_\la^2 )
\xi^{-2 \la t}
\end{eqnarray*}
respectively.
Comparing these two expressions, and using the uniqueness of the Galois
ring transform coefficients, we find
\beql{E.A3}
2( e_\la^2 + f_\la^2 ) = 2 \sum_{0 \le \la_1 < \la_2 \le n-1 \atop
\la_1 + \la_2 = 2 \la} e_{\la_1} e_{\la_2} ~.
\eeq
Now $\mu ( e_\la ) = \tilde{a} ( \la )$,
$\mu (e_{2 \la} ) = \tilde{a} (2 \la ) = \tilde{a} ( \la )^2 = \mu (e_\la )^2$, so \eqn{E.A3} implies
$$\mu (f_\la^2) =
\sum_{0 \le \la_1 \le \la_2 \le n-1 \atop \la_1 + \la_2 = 2 \la}
\tilde{a} ( \la_1 ) \tilde{a} ( \la_2 ) ~,$$
an equation in $GF(2^m)$.
Taking the square root of both sides we obtain
\beql{E.A4}
\mu ( f_\la ) = \sH ( \tilde{a} , \la ) ~.
\eeq
From \eqn{E.A1}, \eqn{E.A2}, \eqn{E.A4} we see that
$$a_\infty + 2 b_\infty + e_0 + 2f_0 + 2 \widehat{b} (0) =0 ~,$$
which implies
$$
a_\infty + \mu (e_0) = a_\infty + \tilde{a} (0) =0 ~,$$

$$
\mu \left( \sqrt{a_\infty e_0} \right) + b_\infty +
\sH ( \tilde{a} , 0 ) + \tilde{b} (0) =0
$$
and the first and third equations of \eqn{E115.0}
now follow.
The second and fourth equations follow easily from the second
equation of \eqn{E.A1}. \hfill $\blacksquare$

For comparison with \eqn{E115.0}, a transform
characterization of Preparata's original code (of the same length
$2^{m+1}$) can be readily derived from the description given by Baker,
van~Lint and Wilson \cite{BLW83}:
a vector $(b, a+b)$ is in this code if and only if
\bary
\label{E115.5}
\tilde{a} (0) + a_\infty & = & 0 \nonumber \\
\tilde{a} (1) & = & 0 \nonumber \\
\tilde{b} (0) + b_\infty & = & 0 \nonumber \\
\tilde{b} (1)^3 & = & \tilde{a} (3) ~.
\eary
The similarity between \eqn{E115.0} and \eqn{E115.5} is evident.
At length 16 (the case $m=3$) the two descriptions must
coincide, since the \NR code is unique (see \S4.5).
This may be verified directly as follows.
\begin{theo}
\label{NR2}
When $m=3$ the `Preparata' code $P$ coincides with Preparata's
original code.
\end{theo}

\PF
It is enough to show that
$\sH ( \tilde{a} , 0) = a_\infty \tilde{a} (0)$ and
$\sH ( \tilde{a} ,1) = \tilde{a} (3)^{1/3} = \tilde{a} (3)^5$.
The cyclotomic cosets mod~7 are
$\{0\}$, $\{1,2,4\}$ and $\{3,5,6\}$.
Hence
$$
\ta (2) = \ta (1)^2 , ~
\ta (4) = \ta (1)^4 , ~
\ta (6) = \ta (3)^2 , ~ \ta (5) = \ta (3)^4 ~.$$
Since $\ta (1) =0$ and $\ta (0) = a_\infty$ are given, we
have
\bara
\sH ( \ta , 0) & = &
\ta (0)^2 + \ta (1) \ta (6) +
\ta (2) \ta (5) + \ta (3) \ta (4) \\
& = & \ta (0)^2 ~=~ \ta (0) a_\infty ~, \\
\sH ( \ta ,1) & = & \ta (0) \ta (1) + \ta (2) \ta (6) + \ta (3) \ta (5) + \ta (4)^2 \\
& = & \ta (3) \ta (5) ~=~ \ta (3)^5 ~,
\eara
as required. \hfill $\blacksquare$

As we have already seen in \S4.5, the appropriate quaternary code in the case $m=3$ is the self-dual octacode.
\subsection{Decoding the quaternary `Preparata' code in the $\ZZ_4$ domain.~~}
There is a very simple decoding algorithm for the
`Preparata' code $\sP$, obtained by working in the $\ZZ_4$ domain.
This is an optimal syndrome decoder:
it corrects all error patterns of Lee weight at most 2,
detects all errors of Lee weight 3, and detects some errors of Lee
weight 4.
A decision tree for the algorithm is shown in Fig.~2.
We use the parity check matrix $H$ given in \eqn{E3.12}, and assume $m$ is odd
and $\ge 3$.

Let $v= (v_\infty , v_0 , \ldots , v_{n-1} ) \in \ZZ_4^{n+1}$ be the received
vector.
The syndrome $Hv'$ has two components, which we write as
\bara
t & = & \sum_{j=0}^{n-1} v_j + v_\infty ~, \\
A+2B & = & \sum_{j=0}^{n-1} v_j \xi^j ~,
\eara
where $A,B \in \sT$.
\begin{figure}
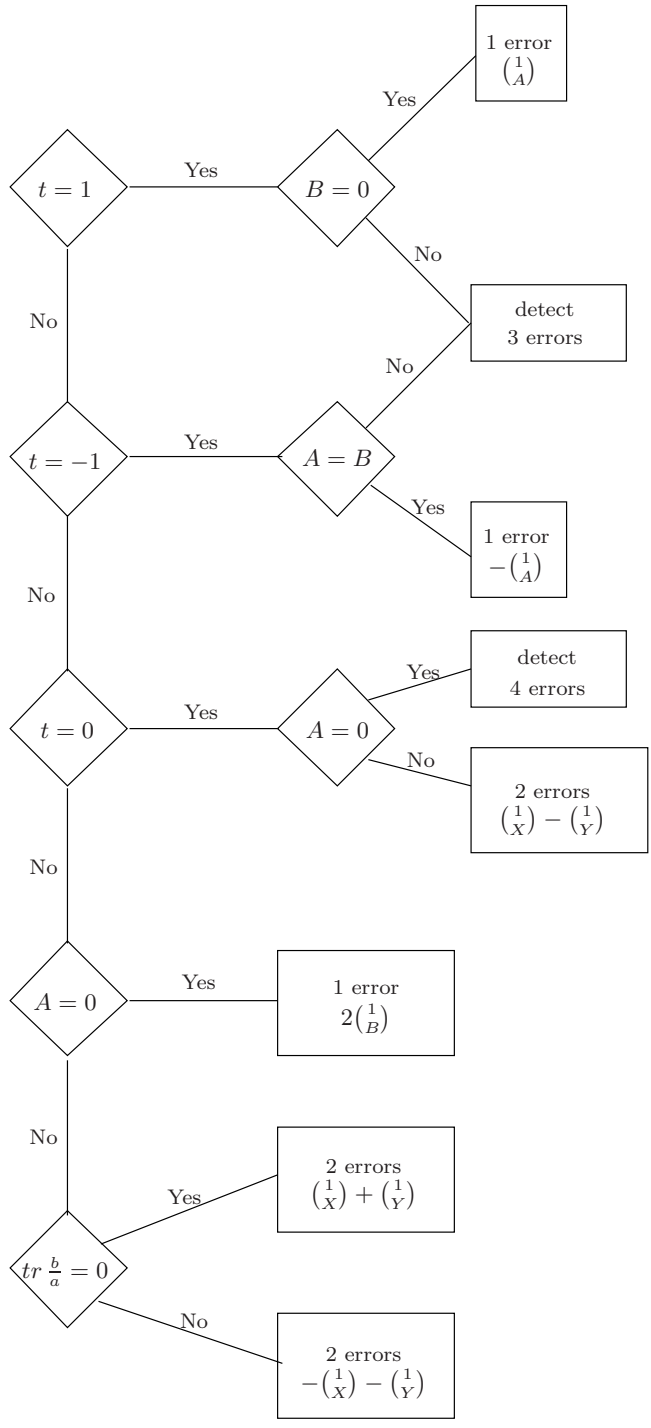

\begin{center}
\input SD3.pstex_t
\end{center}
\caption{Decoding algorithm for `Preparata' code}
\end{figure}

In Theorem~\ref{KW1} we saw that exactly four nonzero weights occur in the Lee
weight distribution of the quaternary \Ker code
$\sK = \sP^\perp$,
and hence also in the Hamming weight distribution of $K$.
It follows that the covering radius of $P$ is at most 4
(\cite{De73};
\cite{MS77},
Theorem~21 of Chap.~6), i.e. the Lee
distance $d_L ( v, \sP )$ from a vector
$v \in \ZZ_4^{n+1}$ to $\sP$ satisfies $d_L (v, \sP ) \le 4$.
Note that $t= \pm 1$ if and only if $d_L (v, \sP ) =1$ or 3.

\vsp
\noindent
{\bf Single errors of Lee weight 1 or 2.}~
If $t=1$ and $B=0$, or if $t= -1$ and $A=B$,
we decide that
there is a single error of Lee weight 1 in column $(1, A)'$.
If $t=1$ and $B \neq 0$, or if $t=-1$ and $A \neq B$, then
$d_L (v, \sP ) =3$.
If $t=2$ and $A=0$,
we decide that there is a single error of Lee weight 2 in column $(1,B)'$.

\vsp
\noindent{\bf Double errors of Lee weight 2.}~
We begin by supposing that $t=0$ and
$$A+2B = X-Y ~,$$
where $X,Y \in \sT$ and $X \neq Y$.
Note that $A \neq 0$ since by (P1) $X-Y$ is invertible.
We have
\bara
A & = & X+Y + 2X^{2^{m-1}} Y^{2^{m-1}} ~, \\
B & \equiv & Y + X^{2^{m-1}} Y^{2^{m-1}} ~~( \bmod~2) ~.
\eara
Let $x,y, a,b$ respectively be the images of $X,Y,A,B$ in $GF(2^m)$ after reduction mod~2 using the map $\mu$.
Then
$$
a = x+y ~, ~~~
b = y + x^{2^{m-1}} y^{2^{m-1}} ~,
$$
which we rewrite as
$$
a = x+y ~, ~~~
(b+y)^2 = xy ~.
$$
The unique solution to these equations is $y=b^2 /a$, $x= a+b^2 /a$.
Note that when $b=0$ or $b=a$, the double error
involves the first column of $H$.
Next we suppose that $t=2$ and that
$$A+2B = X+Y ~,$$
where $X,Y \in \sT$, $X \neq Y$,
$A \neq 0$.
Proceeding as above we find
$$
a = x+y ~, ~~~
b^2 = xy ~,
$$
and so $x$ and $y$ are distinct roots of the
equation
$$u^2 + au + b^2 =0 ~.$$
A necessary and sufficient condition for this equation to have distinct
roots is that
$$tr (b^2 / a^2 ) = tr (b/a) = 0$$
(see \cite{MS77}, Chap.~9, Theorem~15; \cite{LN83}).

Finally we suppose that $t=2$ and
$$A+2B = -X -Y~,$$
where $X,Y \in \sT$, $X \neq Y$,
$A \neq 0$.
We now find that
$$
a = x+y ~,~~~
(b+a)^2 = xy ~,
$$
and so $x$ and $y$ are distinct roots of the 
equation
$$u^2 + au + (a^2 +b^2) =0 ~.$$
A necessary and sufficient condition for this equation to
have distinct roots is that
$$tr \left( \frac{a^2 + b^2}{a^2} \right) = tr
\left( 1 + \frac{b}{a} \right) = 1 + tr \left( \frac{b}{a} \right) =0 ~.
$$
\subsection{Quaternary Reed-Muller codes}
\hsp
In \S2.7 we defined a quaternary code $ZRM (r, m-1)$ whose image under the Gray
map $\phi$ is the binary Reed-Muller code $RM (r,m)$, provided $r \in \{0,1,2, m-1, m \}$.
In this section we define another quaternary Reed-Muller code, $QRM (r,m)$, whose image
under the map $\af$ is $RM(r,m)$ for all $r$, and which includes the \Ker and `Preparata' codes as special cases.

\vsp
\noindent{\bf Definition}.
Let $QRM (0,m)$ be the quaternary repetition code of length $n=2^m$, and for $1 \le r \le m$ let $QRM (r,m)$ be generated by $QRM (0,m)$ together with all vectors of the form
$$(0,T( \la_j ) , T( \la_j \xi^j ) , T( \la_j \xi^{2j} ) , \ldots , T( \la_j \xi^{(n-1)j} ))$$
where $j$ ranges over all representatives of cyclotomic cosets $\bmod~2^m -1$ for which $wt (j) \le r$, and $\la_j$ ranges over $R$.
Then $QRM (r,m)$ is a quaternary code of length $n=2^m$ and type $4^k$, where
$$k = 1+ {\binom{m}{1}} + \cdots + {\binom{m}{r}} ~.$$

\begin{theo}
\label{QRM}
\begin{eqnarray}
\label{E.RM.1}
QRM (1,m) & = & \sK ~, \\
\label{E.RM.2}
QRM (m-2,m) & = & \sP ~, \\
\label{E.RM.3}
\af (QRM (r,m)) & = & RM (r,m) ~, \\
\label{E.RM.4}
QRM (r,m)^\perp & = & QRM (m-r-1,m) ~.
\end{eqnarray}
\end{theo}

\PF
\eqn{E.RM.1} follows from Theorem~\ref{KT}(b), \eqn{E.RM.2} from
\eqn{E.RM.4}, and \eqn{E.RM.3} from \cite[Chap.~13, \S5]{MS77}.
It remains to prove \eqn{E.RM.4}.
This follows from the transform domain characterization of $QRM (r,m)$ as the set of vectors $a$ for which $\widehat{a} ( \la ) =0$ whenever
$wt ( \la ) \le m-1-r$, and $QRM (r,m)^\perp$ as the set of vectors for which $\widehat{a} ( \la ) =0$ whenever $wt ( \la ) \le r$.
(Equivalently, we consider the cyclic codes obtained by deleting the first coordinate, and use the fact that the zeros of a code are the reciprocals of the nonzeros of the dual code.) \hfill $\blacksquare$
\subsection{Automorphism groups.~~}
Consider any system $\Omega$ of linear equations over $\ZZ_4$,
in the variables $c_x$, $x \in \sT$ (see \eqn{E3.2a})), that includes
\bary
\label{ER.1}
\sum_{x \in \sT} c_x & = & 0 ~, \\
\label{ER.2}
\sum_{x \in \sT} c_x x & = & 0 ~,
\eary
together with equations of the form
\beql{ER.3}
2 \left( \sum_{x \in \sT} c_x x^{2^j +1} \right) =0 ~.
\eeq

\begin{theo}
\label{NOT}
The linear system $\Omega$ is invariant under the doubly
transitive group $G$ of `affine' permutations of the form
$$
x \to (ax+b)^{2^m} ~,$$
where $a,b \in \sT$ and $a \neq 0$.
The order of $G$ is $2^m (2^m -1)$.
\end{theo}

\PF
Repeated application of the \For automorphism \eqn{E3.9} to
Eq.~\eqn{ER.2} gives
\beql{ER.4}
\sum_{x \in \sT} c_x x^{2^j} = 0 ~,
\eeq
for all $j$.
It now follows from \eqn{ER.1}, \eqn{ER.2} and \eqn{ER.4}
that
\bara
\sum_{x \in \sT} c_x (ax+b)^{2^m} & = &
\sum_{x \in \sT} c_x (a^{2^m} x^{2^m} +
b^{2^m} + 2a^{2^{m-1}} b^{2^{m-1}} x^{2^{m-1}} ) \\
& = &
a \sum_{x \in \sT} c_x x + b \sum_{x \in \sT} c_x +
2 a^{2^{m-1}} b^{2^{m-1}} \sum_{x \in \sT}
x^{2^{m-1}} ~=~ 0~.
\eara
Finally
\bara
2 \sum_{x \in \sT}
c_x [(ax+b)^{2^m} ]^{2^j +1} & = &
2 \sum_{x \in \sT}
c_x (a^{2^m} x^{2^m} + b^{2^m} )^{2^j +1} \\
& = & 2 \sum_{x \in \sT}
c_x (ax+b)^{2^j +1} \\
& = & 2
\sum_{x \in \sT}
c_x ( a^{2^j +1} x^{2^j +1} +
b^{2^j +1} + a^{2^j} bx^{2^j} + b^{2^j} a x)~=~ 0~.
\eara
Hence $\Omega$ is invariant under $G$. \hfill $\blacksquare$

\vsp
\noindent{\bf Corollary.}
{\em Our quaternary Kerdock, `Preparata', `Goethals', Delsarte-Goethals and
`Goethals-Delsarte' codes are invariant under a 
doubly
transitive group of order $2^{m+1} (2^m -1) m$ generated by
$G$, negation, and the Frobenius map \eqn{E3.9} acting on $\sT$.}

\vsp
\PF
The presence of negation follows from $\ZZ_4$-linearity,
the action of $G$ from Theorem~\ref{NOT} and the Frobenius map from Eq.~\eqn{ER.4}. \hfill $\blacksquare$

\vsp
\noindent{\bf Remarks.}
By the automorphism group $\mbox{Aut} (C)$ of a binary nonlinear code $C$
we will mean the set of all coordinate permutations that preserve the code.
It is easy to see that if $C = \phi ( \sC )$ is the binary image of a linear
quaternary code $\sC$, then $\mbox{Aut} ( \sC )$ is isomorphic to a subgroup of $\mbox{Aut} (C)$.

The automorphism groups of the binary Nordstrom-Robinson, Kerdock,
classical Preparata, and Delsarte-Goethals codes are
known (Berlekamp \cite{Ber71}, Carlet \cite{Car91}, \cite{Car92}, \cite{Car93},
Kantor \cite{Ka82a,Ka83}).
For odd $m \ge 5$ these groups have the same orders as those in the Corollary.

We conclude that, {\em for odd $m \ge 5$, the groups mentioned in the Corollary are the full automorphism groups of these quaternary codes.}
(For the `Preparata' codes we use the fact that they
have the same automorphism group as their duals.)

The case $m=3$ is exceptional.
The quaternary octacode has an automorphism group of order 1344 (Conway and Sloane \cite{CS93}), whereas the group of
the binary Nordstrom-Robinson code
has order 80640 (Berlekamp \cite{Ber71}, see also Conway and Sloane \cite{CS90}).
\subsection{A new family of distance regular graphs of diameter 4.~~}
As before, $P = \phi ( \sP ) = \phi ( \sK^\perp )$ denotes our `Preparata'
code of length $N= 2^{m+1}$, with $m$ odd $\ge 3$.

\vsp
\noindent{\bf Definition.}
A {\em $\ZZ_4$-coset} of $P$ is the image under $\phi$ of a coset of $\sP$ in
$\ZZ_4^{N/2}$.
We construct a graph $\Gamma_m$ on the $\ZZ_4$-cosets of $P$ by
joining two cosets by an edge if they are the images of cosets
$x+ \sP$, $y + \sP$ such that $x-y+ \sP$ has minimal Lee weight 1.

Let $\Pi$ denote the partition of $\ZZ_2^N$ into $\ZZ_4$-cosets of $P$.
Then $\Gamma_n$ can be thought of as the quotient graph
(\cite[\S11.1.B]{BCN89})
of the $N$-hypercube by the partition $\Pi$.

The aim of this section is to show that $\Gamma_m$ is distance regular
and to compute its distance distribution
diagram and eigenmatrix {\bf P}.
For this purpose we need certain regularity properties of $P$ and $\Pi$.

If $C$ is a binary code of length $N$, its outer distribution matrix
$B = (B_{x,j} )$ is the $2^N \times (N+1)$ matrix
with typical entry
$$B_{x,j} = \left|
\{ y \in C: ~d(x,y) = j \} \right|$$
(Delsarte \cite{De73}).
In other words the rows of $B$ are the weight distributions
of the translates of $C$.

A {\em code} $C$ of covering radius $r$ is said to be
{\em completely regular}
\cite{De73a} if $B$ contains exactly $r+1$ distinct rows.
A {\em partition} $\Pi$ of $\ZZ_2^N$ into cosets is said to be
{\em completely regular} if all members of the partition are
completely regular with the same matrix.
\begin{lem}
\label{A1}
The covering radius of $P$ is 4.
\end{lem}

\PF
In the previous section we saw that it is at most 4.
But $P$ is contained in a code with the same weight distribution as an extended
Hamming code (Theorem~\ref{PD}), and so by the supercode
lemma \cite{CKMS} the
covering radius is at least 4. \hfill $\blacksquare$
\begin{lem}
\label{A2}
The codewords of weight 6 in $P$ form a $3- (2^{m+1} , 6, (2^{m+1} -4)/3)$ design.
\end{lem}

\PF
The proof of Theorem~33 of \cite{MS77}, Chap.~15 can be used, since it
depends only on the annihilator polynomial of $P$. \hfill $\blacksquare$

\begin{theo}
\label{A}
The `Preparata' code $P$ is completely regular.
\end{theo}

\PF
The well-known recurrence relation between the columns of $B$
(\cite{De73}, \cite{MS77}) has order 4, by Lemma~\ref{A1},
and so it is sufficient to check that $B_{x,j}$ can take at most five different values for fixed
$x \in \ZZ_2^N$ and $0 \le j \le 4$.
If $d(x,P) \le 2$, the fact that $P$ has minimal distance 6 shows that
$B_{x,j}$ is either 0 or 1.
If $d(x,P) = 3$ then Lemma~\ref{A2} shows that $B_{x,3} = (N-1)/3$.
Clearly $B_{x,0} = B_{x,1} = B_{x,2} =0$.
Finally if $d(x,P) =4$ then $B_{x,0} = B_{x,1} = B_{x,2} = B_{x,3} =0$. \hfill $\blacksquare$

As in \cite{CM91} it will be noticed that $P$ is neither linear, perfect,
nor uniformly packed, and so (in the notation of Levenshtein \cite{Lev})
is not a design of
Delsarte type (i.e. $d \ge 2s' -1$);
$P$ is a highly nontrivial example of a completely regular code.
Furthermore the $\ZZ_4$-linearity of $P$ and the properties
of $\phi$ show that each $\ZZ_4$-coset of $P$ is completely
regular with the same outer distribution matrix.
Hence $\Pi$ is completely regular.
The next result follows immediately from Theorems~11.1.6 and
11.1.5 of \cite{BCN89}.
\begin{theo}
\label{B}
The graph $\Gamma_m$ is distance regular on $N^2$ vertices with
diameter 4 and degree $N$.
\end{theo}

We now proceed to a more detailed study of the parameters of
$\Gamma_m$.
Recall that the valencies $v_j$ are the numbers of points at distance $j$
from a given point.
The intersection numbers $a_j$, $b_j$, $c_j$ are defined in
Chapter~4 of \cite{BCN89}.
\begin{lem}
\label{C1}
$\Gamma_m$ is bipartite.
\end{lem}

\PF
Let us take a parity check matrix $H$ of the form \eqn{E3.12}
for $\sP$.
For a coset $x+ \sP$ let $Hx'$ be the associated syndrome
and let $\nu (x)$ be the leading bit of $Hx'$.
Then $\nu$ is a map from the vertices of $\Gamma_m$ onto
$\{0,1\}$.
Let $X_j$ be the set $\nu^{-1} (j)$, $j=0,1$.
Since $\nu (x) =1$ if $x$ has
weight 1, two cosets with the same image under $\nu$ cannot have
adjacent images in $\Gamma_m$. \hfill $\blacksquare$
\begin{lem}
\label{C2}
If $x \in \ZZ_2^n$ is at distance 4 from $P$, then
$B_{x,4} = N(N-1)/12$.
\end{lem}

\PF
From \cite{So90} and the fact that $P$ has size $2^N /N^2$ and four
nonzero dual distances $d'_1$, $d'_2$, $d'_3$, $d'_4$ we obtain
$$B_{x,4} = \frac{2^4}{4!N^2} \prod_{j=1}^4 d'_j ~.$$
The desired result then follows from
$d'_1 = (N- \sqrt{N})/2$,
$d'_2 = N/2$, $d'_3 = (N+\sqrt{N})/2$, $d'_4 =N$.~~\hfill $\blacksquare$
\begin{theo}
\label{C}
The valencies of $\Gamma_m$ are $v_0 =1$, $v_1 =N$, $v_2 = {\binom{N}{2}}$,
$v_3 = \frac{N(N-2)}{2}$, $v_4 = \frac{N-2}{2}$.
The intersection numbers of $\Gamma_m$ are $b_0 =N$, $c_1 =1$,
$b_1 =N-1$, $c_2 =2$,
$b_2 = N-2$, $c_3 = N-1$,
$b_3 =1$, $c_4 =N$.
Furthermore $a_j =0$ for $j=0,1,2,3,4$.
\end{theo}

\PF
By Lemma~\ref{C1},
$\Gamma_m$ is bipartite, hence without circuits of odd length.
Therefore $a_j =0$ for $0 \le j \le 4$.

The intersection numbers add up to the degree, so
$N= b_j + c_j$ for $0 \le j \le 4$,
and it only remains to calculate the $c_j$.
The values of $c_1$ and $c_2$ are clear from the double-error-correcting
character of $\sP$.
Finally $c_3$ and $c_4$ are computed from the formula of Theorem~11.1.8 of
\cite{BCN89} by observing that $e_{l,j} = B_{x,j}$ if
$d(x,P) =l$.
Moreover $e_{3,3} = (N-1)/3$ by Lemma~\ref{A2} and $e_{4,4} = N(N-1)/12$
by Lemma~\ref{C2}.
The intersection numbers of the $N$-cube are
well known to be $a_j =0$, $b_j = N-j$,
$c_j = j$. \hfill $\blacksquare$

\vsp
\noindent{\bf Corollary}.
{\em
The eigenmatrix {\bf P} for $\Gamma_m$ is
}
$$
P = \pmatrix{
1 & N & {\dis\binom{N}{2}} & \df{N(N-2)}{2} & \df{N-2}{2} \cr
1 & \sqrt{N} & 0 & - \sqrt{N} & -1 \cr
1 & 0 & - \df{N}{2} & 0 & 1 - \df{N}{2} \cr
1 & - \sqrt{N} & 0 & \sqrt{N} & -1 \cr
1 & -N & {\dis\binom{N}{2}} & - \df{N(N-2)}{2} & \df{N}{2} -1 \cr}
$$

\vsp
\PF
See \cite{BCN89}, \S4.1.B, or \cite{CM91},
Proposition~3.17. \hfill $\blacksquare$

\vsp
\noindent{\bf Remarks.}
1)~~Let $R_j$ denote the $j$th class of the association
scheme corresponding to $\Gamma_m$.
Then $R_1 + R_3$ has only three eigenvalues,
$N^2$, 0, $-N^2$, and is a strongly regular graph
isomorphic to the complete bipartite graph
$K_{N^2 /2 , \, N^2 /2}$.
It would be interesting to see if $R_3$ is also distance regular.

\noindent
2)~$\Pi$ is a 4-partition design in the sense of \cite{CCD92}, \cite{CM91}.

\section{Goethals, Delsarte-Goethals, and other codes}
\hsp
It is natural to wonder how the constructions of $\sK$ and $\sP$
can be generalized.
We have already seen one generalization in \S5.4.
Another generalization is
to replace \eqn{E3.12} by the matrix
\beql{E6.1}
\left[
\begin{array}{cccccc}
1 & 1 & 1 & 1 & \cdots & 1 \\
~~ \\ [-.09in]
0 & 1 & \xi & \xi^2 & \cdots & \xi^{(n-1)} \\
~~ \\ [-.09in]
0 & 2 & 2 \xi^3 & 2 \xi^6 & \cdots & 2 \xi^{3(n-1)} \\
~~ \\ [-.09in]
~ & \cdot & \cdot & \cdot & \cdots & ~ \\
~~ \\ [-.09in]
0 & 2 & 2 \xi^{1+2^j} & 2 \xi^{2(1+2^j )} & \cdots & 2 \xi^{(1+2^j) (n-1)} \\
~~ \\ [-.09in]
~ & \cdot & \cdot & \cdot & \cdots & ~ \\
~~ \\ [-.09in]
0 & 2 & 2 \xi^{1+2^r} & 2 \xi^{2(1+2^r )} & \cdots & 2 \xi^{(1+2^r )(n-1)}
\end{array}
\right] ~,
\eeq
where $1 \le r \le (m-1)/2$.
Again we assume $m$ is odd.
\begin{theo}
\label{DG}
(a)~The quaternary code of length $2^m$ with generator matrix \eqn{E6.1}
has type $4^{m+1} 2^{rm}$ and minimal Lee weight
$2^m - 2^{m- \delta}$, where
$\delta = \frac{m+1}{2} -r$.
The binary image under the Gray map \eqn{E1} is the \DG code
$DG (m+1 , \delta )$
(\cite{DG75}; \cite{MS77}, Chap.~15).
(b)~The dual code, with parity check matrix \eqn{E6.1}, has a binary image
with minimal distance $8$ and the same weight distribution as the Goethals-Delsarte code
$GD (m+1 , r+2)$ 
defined by Hergert \cite{He90}.
In particular, for $r=1$ this produces a binary code
$G$
with
the same weight
distribution as the Goethals code $\sI (m+1)$
(\cite{Go74}; \cite{Go76}; \cite{MS77}, Chap.~15).
\end{theo}

\PF
(a)~Comparing Eqs.~(37) and (34) of \cite{MS77}, Chap.~15,
we see that the difference between the Kerdock code and the \DG code comes
from the vectors $(c,c)$, where $c$ belongs to the code defined
by Eq.~(31) of that chapter.
We already know from Theorem~\ref{K} that the first two rows of \eqn{E6.1}
produce the Kerdock code, and it is easily seen that the remaining
rows produce the required $(c,c)$ words.
(b)~This follows because the Goethals-Delsarte code is by construction (see Hergert \cite{He90}) a distance invariant codes whose weight enumerator is the
MacWilliams transform of the \DG code.
The minimal Lee distance of these dual codes is no more than 8,
since they contain words of shape $2^4$, corresponding to the doubles
of words in the extended Hamming code defined by the binary
images of the first two rows of \eqn{E6.1}.
That the minimal Lee distance is at least 8 follows from
Theorem~\ref{G8} below. \hfill $\blacksquare$

\vsp
\noindent{\bf Remarks.}
1)~There are also transform-domain characterizations of some of these codes.
For the `Goethals' codes and the dual codes defined in part (b) of Theorem~\ref{DG}, add to \eqn{E115.0}
the conditions
$$\tilde{a} (1+2^i) =0 , ~~~ i=1,2, \ldots , r ~.$$
For the original Goethals codes \cite[p.~477]{MS77}, replace \eqn{E115.5} by
\begin{eqnarray*}
\tilde{a} (0) & + & a_\infty ~=~ 0 ~, \\
& ~ & \tilde{a} (1) ~=~ 0 ~, \\
\tilde{b} (0) & + & b_\infty ~=~ 0 ~, \\
\tilde{a} (r) & = & \tilde{b} (1)^r ~, \\
\tilde{a} (s) & = & \tilde{b} (1)^s ~,
\end{eqnarray*}
where $r = 1+2^{t-1}$, $s=1+2^t$, and $a$, $b$ are binary vectors of length $n=2^{2t+1}$.

\noindent
2)~For the automorphism groups of these codes, see Section~5.5.

\noindent
3)~Our `Goethals' code is thus defined as $G = \phi ( \sG )$, where $\sG$
is the quaternary code with parity check matrix
\beql{E200}
\left[
\begin{array}{cccccc}
1 & 1 & 1 & 1 & \cdots & 1 \\
~~ \\ [-.09in]
0 & 1 & \xi & \xi^2 & \cdots & \xi^{(n-1)} \\
~~ \\ [-.09in]
0 & 2 & 2 \xi^3 & 2 \xi^6 & \cdots & 2 \xi^{3(n-1)}
\end{array}
\right] ~.
\eeq
We end by giving a direct proof that this code has minimal distance 8.
\begin{theo}
\label{G8}
The minimal distance of the `Goethals' code $G = \phi ( \sG )$ of
length $2^{m+1}$, $m$ odd $\ge 3$, is 8.
\end{theo}

\PF
Since $\sG \subseteq \sP$, the minimal distance $d$ is at least 6.
Suppose, seeking a contradiction, that $c= (c_\infty, c_0 , \ldots , c_{n-1} )$ is
a codeword of type $( \pm 1 )^{n_1} 2^{n_2} 0^{n+1-n_1 -n_2}$,
where $n= 2^m -1$, $n_1 + 2n_2 =6$.
Write $c= 2c_0 + c_1$,
where $2c_0$ is a vector of type $2^{n_2} 0^{n+1-n_2}$ and $c_1$
is a vector of type $( \pm 1 )^{n_1} 0^{n+1-n_1}$.
Then $c_1$ is orthogonal to every row of the matrix
$$
2 \left[
\matrix{
1 & 1 & 1 & \cdots & 1 \cr
0 & 1 & \xi & \cdots & \xi^{n-1} \cr
0 & 1 & \xi^3 & \cdots & \xi^{3(n-1)} \cr}
\right ] ~,
$$
and so $\af (c_1)$ is in the extended
double-error-correcting BCH code of length $2^m$.
It follows that $n_2 =0$ and $n_1 =6$,
and in fact that $c$ must be of the type
$1^3 3^3 0^{n-5}$ or $\pm (1^5 3~ 0^{n-5} )$.

\vsp
\noindent
{\bf Case~1}.
$c$ is of type $1^3 3^3 0^{n-5}$. The
automorphism group of $\sG$ is doubly
transitive on the coordinate positions (Theorem~\ref{NOT}),
so we may assume $c_\infty = -1$.
Thus $c$ determines a solution to the equations
\bara
X_1 + X_2 + X_3 & = & Z_1 + Z_2 ~, \\
X_1^3 + X_2^3 + X_3^3 & \equiv & Z_1^3 + Z_2^3 ~~~ (\bmod~2) ~,
\eara
where $X_1$, $X_2$, $X_3$, $Z_1$, $Z_2$ are distinct nonzero
elements of $\sT$.
If $x_1$, $x_2$, etc., are the images of these elements in $GF(2^m)$ under $\mu$, we have
\bara
x_1 + x_2 + x_3 & = & z_1 + z_2 ~, \\
x_1^3 + x_2^3 + x_3^3 & = & z_1^3 + z_2^3 ~, \\
x_1 x_2 + x_2 x_3 + x_3 x_1 & = & z_1 z_2 ~.
\eara
But this implies
\bara
x_1 x_2 x_3 & = &
(x_1 + x_2 + x_3 )^3 +
(x_1^3 + x_2^3 +x_3^3 ) \\
& ~ & ~~~+ ~ (x_1 + x_2 + x_3 )
(x_1 x_2 + x_2 x_3 + x_3 x_1 ) \\
& = &
(z_1 + z_2 )^3 + z_1^3 +z_2^3 +
z_1 z_2 ( z_1 + z_2 ) ~=~ 0 ~,
\eara
which is a contradiction.

\vsp
\noindent{\bf Case~2}.
$c$ is of type $1^5 3~0^{n-5}$.
By using the automorphism group we may suppose $c_\infty =3$, $c_0 =1$.
Thus $c$ determines a solution to
\bara
X_1 + X_2 + X_3 & = & -1 - Z_1 ~, \\
X_1^3 + X_2^3 + X_3^3 & \equiv & -1-Z_1^3 ~~~(\bmod~2) ~,
\eara
where $X_1$, $X_2$, $X_3$, 1, $Z_1$ are distinct nonzero elements of $\sT$.
Proceeding as before we find
\bara
x_1 + x_2 + x_3 & = & 1 + z_1 ~, \\
x_1^3 + x_2^3 + x_3^3 & = & 1+ z_1^3 ~, \\
x_1 x_2 + x_2 x_3 + x_3 x_1 & = & 1 + z_1 + z_1^2 ~.
\eara
For $i=1,2,3$ let $y_i = x_i +1+z_1$.
This change of variables produces the equations
\bary
\label{EP5.1}
y_1 + y_2 + y_3 & = & 0 ~, \nonumber \\
y_1^3 + y_2^3 + y_3^3 & = & z_1 (1+z_1) ~, \\
y_1 y_2 + y_2 y_3 + y_3 y_1 & = & z_1 ~. \nonumber
\eary
Now write $y_2 = ay_1$,
$y_3 = (1+a) y_1$, so that
\bara
y_1^3 (a+a^2) & = & z_1 + z_1^2 ~, \\
y_1^2 (1+a+a^2) & = & z_1 ~, \\
\eara
and also $y_1 \neq 0$.
It follows that
$$y_1^2 (1+ a^2 + a^4) + y_1 (a+a^2)
+ (1+a+a^2) =0 ~.$$
Setting $s = a+a^2$, we obtain the quadratic equation
\beql{EP5.2}
s^2 + \frac{(1+y_1)}{y_1^2} s +
\frac{1+y_1^2}{y_1^2} =0 ~.
\eeq
We shall prove that this equation has no solution.
First observe that $y_1 \neq 1$, so the equation does not have a
double root.
Suppose the equation has two distinct roots.
It follows from \eqn{EP5.1} that there exist distinct
nonzero elements $Y_1$, $Y_2$,
$Y_3$, $Y'_2$, $Y'_3$, $Z_1^{1/2}$ of
$\sT$ such that
$$Y_1 + Y_2 + Y_3 = 2 Z_1^{1/2} = Y_1 + Y'_2 + Y'_3 ~.$$
However, this implies the existence of codewords in the `Preparata' code of
type $1^2 3^2~0^{n-3}$, which is not the case.
Hence \eqn{EP5.2} has no solutions and the proof is complete. \hfill $\blacksquare$

We are presently investigating other generalizations of \eqn{E3.12}.
\section{Conclusions}
\hsp
The classical theory of cyclic codes, which includes BCH,
Reed-Solomon, Reed-Muller codes, etc., regards these codes as
ideals in polynomial rings over finite fields.
Some famous nonlinear codes found by Nordstrom-Robinson, Kerdock,
Preparata, Goethals and others, more powerful than any linear
codes, cannot be handled by this machinery.
We have shown that when suitably defined all these codes are ideals in polynomial rings over
the ring of integers mod~4.
This new point of view should completely transform the study of
cyclic codes.
\section*{Acknowledgements}
\hsp
We thank Claude Carlet, Pascale Charpin, Dave Forney,
Tor Helleseth,
Vladimir Levenshtein,
Kyeongcheol Yang and Viktor Zinoviev for helpful discussions
and comments.

\clearpage
\begin{center}
{\bf Key~words and phrases}
\end{center}
Kerdock code, Preparata code,
Nordstrom-Robinson code, Goethals code,
Delsarte-Goethals code,
Goethals-Delsarte code,
octacode, nonlinear codes, quaternary codes,
Reed-Muller codes, cyclic codes,
completely regular codes.
\clearpage
\begin{center}
\section*{List of Figure Captions}
\end{center}
Figure~1.
Gray encoding of quaternary symbols and QPSK phases.

\vspace{.1in}
\noindent
Figure~2.
Decoding algorithm for `Preparata' code.

\end{document}